\documentclass[a4paper,10pt]{article}

\usepackage{amsmath}
\usepackage{amssymb,amstext}
\usepackage{amsthm}
\usepackage{algorithm, algorithmic}
\usepackage{url}
\usepackage{array}
\usepackage{booktabs}
\usepackage{psfrag}
\usepackage{graphicx}
\usepackage[footnotesize,nooneline,bf]{caption}
\newcommand{\R}{\mathbb{R}}
\newcommand{\N}{\mathbb{N}}

\newtheorem{theorem}{Theorem}[section]
\newtheorem{lemma}[theorem]{Lemma}
\newtheorem{definition}[theorem]{Definition}
\newtheorem{corollary}[theorem]{Corollary}

\theoremstyle{definition}
\newtheorem{example}[theorem]{Example}
\newtheorem{remark}[theorem]{Remark}

\DeclareMathOperator{\argmin}{argmin}
\DeclareMathOperator{\Span}{span}
\DeclareMathOperator{\rank}{rank}

\DeclareMathOperator{\treeop}{tree}

\newcommand{\supp}[1]{\suppop({#1})}
\newcommand{\suchthat}{\, : \,}
\newcommand{\card}[1]{\#\,{#1}}
\newcommand{\T}{^{\rm {T}}}
\newcommand{\norm}[2]{\lVert{#1}\rVert_{#2}}
\newcommand{\Norm}[2]{\big\lVert{#1}\big\rVert_{#2}}
\newcommand{\abs}[1]{\lvert{#1}\rvert}

\newcommand{\define}{:=} 
\newcommand{\Ref}[1]{\refop(#1)}
\newcommand{\Tree}[1]{\treeop(#1)}

\newcommand{\iprod}[2]{\langle{#1},{#2}\rangle}
\newcommand{\hprod}[2]{({#1},{#2})}
\newcommand{\eprod}[3]{({#1},{#2})_{#3}}

\newcommand{\U}{\mathbb{U}}   
\newcommand{\V}{\mathbb{V}}   
\newcommand{\HS}{\mathbb{H}}  
\newcommand{\Lop}{L}          
\newcommand{\MutInc}[1]{\mathcal{M}({#1})}  

\graphicspath{{Pictures/}}

\title{Sparse Approximate Solution of Partial Differential Equations%
\footnote{Supported by the \textit{Deutsche Forschungsgemeinschaft} through the  DFG
Research Center \textsc{Matheon} \textit{Mathematics for key technologies} in Berlin.}}

\author{Sadegh Jokar%
\footnote{Institut f\"ur Mathematik, TU Berlin, Str.\ des 17.\ Juni 136,
 10623 Berlin, Germany.
 \texttt{$\{$jokar,mehrmann,yserenta$\}$@math.tu-berlin.de}.}
\and Volker Mehrmann\footnotemark[2]
\and Marc Pfetsch%
\footnote{Zuse Institute Berlin, Takustr.\ 7, 14195 Berlin, Germany. \texttt{pfetsch@zib.de}}
\and  Harry Yserentant\footnotemark[2]}
\date{	}

\begin{document}

\maketitle

\begin{abstract}
  \noindent A new concept is introduced for the adaptive finite element
  discretization of partial differential equations that have a sparsely
  representable solution.  Motivated by recent work on compressed sensing,
  a recursive mesh refinement procedure is presented that uses linear
  programming to find a good approximation to the sparse solution on a
  given refinement level.  Then only those parts of the mesh are refined
  that belong to large expansion coefficients.  Error estimates for this
  procedure are refined and the behavior of the procedure is demonstrated
  via some simple elliptic model problems.
\end{abstract}



\section{Introduction}

The sparse representation of functions via a linear combination of a small
number of basic functions has recently received a lot of attention in
several mathematical fields such as approximation theory
\cite{CohDD06a,KunR06,Rau06a,Rau07a} as well as signal and image processing
\cite{Can06,CanR06,CanRT06,CanRT06b,CanT05,CanT06,CheDS99,Don06a,Don06c,Don06b,DonE03,DonET06,DonH01,DonTDS06}.
In terms of representations of functions, we can describe the problem as
follows. Consider a linearly dependent set of $n$ functions $\phi_i$, $i =
1,2,\dots,n$, (\emph{a dictionary}~\cite{Chr03}) and a function $f$
represented as
\[
  f = \sum_{i=1}^n x_i\, \phi_i.
\]
Since the set of functions is not linearly independent, this representation
is not unique and we may want to determine the sparsest representation,
i.e., a representation with a maximal number of vanishing coefficients
among $x_1, \dots, x_n$.
In the setting of numerical linear algebra, this problem can be formulated
as follows. Consider a linear system
\begin{equation}\label{eq:LS}
  \Phi x = b,
\end{equation}
with $\Phi \in \R^{m,n}$, where $m \leq n$ and $b \in \R^m$. The columns of
the matrix $\Phi$ and the right hand side $b$ represent the functions
$\phi_i$ and the function $f$, respectively, with respect to some basis of
the relevant function space. The problem is then to find the sparsest
possible solution $x$, i.e., $x$ has as many zero components as possible.
This optimization problem is in general NP-hard~\cite{GarJ79,Nat95}.
Starting from the work of~\cite{CheDS99}, however, a still growing number
of articles have developed sufficient conditions that guarantee that an
(approximate) sparse solution $\hat{x}$ to~\eqref{eq:LS} can be obtained by
solving the linear program
\[
  \min\ \norm{x}{1},\ s.t.\ \Phi x = b \ (\norm{\Phi x-b}{} \leq \epsilon),
\]
which can be done in polynomial time~\cite{JokP08,Kar84,Kha79}.  We will give a
brief survey of this theory in Section~\ref{subsec:CompressedSensing}.

In the literature, the development has mostly focused on the construction
of appropriate coding matrices $\Phi$ that allow for the sparse
representation of a large class of functions (signals or images).
Furthermore, properties of the columns of the matrix (or the dictionary)
have been investigated, which guarantee that the computation of the sparse
solution can be done efficiently via a linear programming approach, see,
for instance, \cite{CanT06,JokP08}. Often the term \emph{compressed sensing} is
used for this approach.

In this paper we consider a related but different problem.  We are
interested in the numerical solution of partial differential equations
\[
  \Lop u = f,
\]
with a differential operator $\Lop$, to be solved in a domain $\Omega
\subset \R^d$ with smooth boundary $\Gamma$ and appropriate boundary
conditions given on $\Gamma$.

Considering a classical Galerkin or Petrov-Galerkin finite element
approach, see e.g.\ \cite{Bra07}, one seeks a solution $u$ in some function
space $\U$ (which is spanned by $\phi_1, \dots, \phi_n$), represented as
\begin{equation}\label{eq:rep}
  u = \sum_{i=1}^n u_i\, \phi_i.
\end{equation}
Again we are interested in sparse representations with a maximal number of
vanishing coefficients $u_i$. In contrast to the cases discussed before,
here we would like to construct the space $\U$ and the basis functions
$\phi_i$ in the finite element discretization in such a way that first of
all a sparse representation of the solution to~\eqref{eq:rep} exists and
second that it can be determined efficiently. Furthermore, it would be
ideal if the functions $\phi_i$ could be constructed in a multilevel or
adaptive way.

The usual approach to achieve this goal is to use local a posteriori error
estimation to determine where a refinement, i.e., the addition of further
basis functions is necessary. For example, in the dual weighted residual
approach~\cite{BecR01} this is done by solving an optimization problem for
the error.

In this article, we examine the possibility to use similar approaches as
those used in compressed sensing, i.e., to use $\ell_1$-minimization and
linear programming to perform the adaptive refinement in the finite element
method in such a way that the solution is sparsely represented by a linear
combination of basis functions. In order to achieve this goal, we propose
the following framework.

We determine $u \in \U$ as the solution of the weak formulation
\[
(v, \Lop u - f) = 0\ \text{ for all }\ v \in \V.
\]
Here, $\V$ is a space of test functions and $(\cdot{,}\cdot)$ is an
appropriate inner product. In the simplest version of a two-level approach,
we construct finite dimensional spaces of coarse and fine basis
functions 
$\U^n_1 \subset \U^N_1 \subset \U$ and corresponding spaces for coarse and
fine test functions $\V^n_1 \subset \V^N_1 \subset \V$. Then we determine
the sparsest solution in $\U^N_1$, such that
\[
  (v, \Lop u - f) = 0\ \text{ for all } v\in \V_N^1 \setminus \V_n^1
\]
via the solution of an underdetermined system of the form~\eqref{eq:LS}.
Based on the sparse solution, we determine new coarse and fine spaces
$\U^n_2 \subset \U^N_2\subset \U$, $\V^n_2 \subset \V^N_2\subset \V$, and
iterate this procedure.

This framework combines the ideas developed in compressed sensing with
well-known concepts arising in adaptive and multilevel finite element
methods \cite{CohDD01a}. But instead of using local and global error estimates to obtain
error indicators by which the grid refinement is controlled, here the
solution of the $\ell_1$-mini\-miz\-ation problem is used to control the grid
refinement and adaptivity.

Many issues of this approach have, however, not yet been resolved, in
particular, the theoretical analysis of this approach (see
Section~\ref{sec:rippde}). We see the following potential advantages and
disadvantages of this framework. On the positive side, the
$\ell_1$-mini\-miz\-ation approach allows for an easy automation. We will
demonstrate this with some numerical examples in
Section~\ref{sec:numerics}. On the downside, the analysis of the approach
seems to be hard even for classical elliptic problems, see
Section~\ref{sec:rippde} and due to the potentially high complexity of the
linear programming methods this approach will only be successful, if the
procedure needs only a few levels and a small sparse representation of the
solution exists, see Section~\ref{sec:numerics}.

\section{Notation and Preliminaries}
\label{sec:notation}

\subsection{Notation}\label{sec:not}

For $m,n \in \N =\{1,2,\dots\}$, we denote by $\R^{m,n}$ the set of real $m
\times n$ matrices, and by $I_n$ the $n \times n$ identity matrix.
Furthermore, we denote the Euclidean inner product on $\R^n$ by
$\iprod{\cdot}{\cdot}$, i.e., for $x,y \in \R^n$,
\[
\iprod{x}{y} = \sum_{j=1}^n x_j\, y_j.
\]
For $1 \leq p \leq \infty$, the $\ell_p$-norm of $x \in \R^n$ is defined by
\[
  \norm{x}{p} \define \big(\sum_{j=1}^n \abs{x_j}^p\big)^{\frac{1}{p}},
\]
with the special case
\[
\norm{x}{\infty} \define \max_{j \in \{1,\dots,n\}} \abs{x_j},
\]
if $p = \infty$.

The definition of $\norm{\cdot}{p}$ can also be formally extended to the
case that $0 \leq p < 1$. For $0 < p < 1$, $\norm{\cdot}{p}$ is only a
quasi-norm, since the triangle inequality is not satisfied, but still a
\emph{generalized triangle inequality} holds, i.e., for every $x, y \in
\R^n$ one has
\[
\norm{x+y}{p}^p \leq \norm{x}{p}^p + \norm{y}{p}^p.
\]
Finally, for $p=0$ and $x \in \R^n$, we introduce the notation
\[
\norm{x}{0} \define \card{\supp{x}},
\]
where $\supp{x} \define \{j \in \{1, \dots, n\} \suchthat x_j \neq 0\}$ is
the \emph{support} of~$x$. Hence, $\norm{x}{0}$ counts the number of
nonzero entries of~$x$. Note that in this case even the homogeneity is
violated, since for $\alpha \neq 0$ we have $\norm{\alpha x}{0} =
\norm{x}{0}$.

For a symmetric positive definite matrix $A = A\T \in \R^{n,n}$, we
introduce the \emph{energy inner product}
\[
\eprod{u}{v}{A} \define \iprod{u}{Av}
\]
and the induced \emph{energy norm}
\[
\norm{x}{A} \define \sqrt{\eprod{x}{x}{A}}.
\]
Every symmetric positive definite matrix $A \in \R^{n,n}$ has a unique
symmetric positive definite \emph{square root} $B \define A^{\frac{1}{2}}$, with
$A = B^2 = B\T B$ satisfying the relation~\cite{HorJ86}:
\[
\norm{x}{A} = \norm{Bx}{2}.
\]

\subsection{Sparse Representation and Compressed Sensing}
\label{subsec:CompressedSensing}

In this part we survey some recent results on sparse representations of
functions based on the solution of underdetermined linear systems via
$\ell_1$-mini\-mi\-za\-tion. We also discuss the recently introduced
concept of compressed sensing.

\begin{definition}[\cite{CanRT06b}]\label{def:RIP}
  Let $\Phi\in \R^{m,n}$ with $m \leq n$ and $k \in \{1, \dots, n\}$. The
  \emph{$k$-restricted isometry constant} is the smallest number $\delta_k$, such that
  \begin{equation}\label{eq:RIPeq1}
    (1-\delta_k) \norm{x}{2}^2 \leq \norm{\Phi x}{2}^2 \leq (1+\delta_k) \norm{x}{2}^2
  \end{equation}
  for all $x \in \R^n$ with $\norm{x}{0} \leq k$.
\end{definition}

If $\Phi$ in Definition~\ref{def:RIP} is orthonormal, then clearly
$\delta_k = 0$ for all~$k$. Conversely, if the constant $\delta_k$ is close
to $0$ for some matrix~$\Phi$, every set of columns of $\Phi$ of
cardinality less than or equal to $k$ behaves like an orthonormal system.
In the case that $0 \leq \delta_k<\sqrt{2}-1$ for large enough $k$, we say that the matrix $\Phi$ has the
\emph{restricted isometry property}~\cite{Can08,DoVore07a}.

For $\Phi\in \R^{m,n}$ with $m \leq n$, a vector of the form $b=\Phi x$
\emph{represents} (encodes) the vector $x$ in terms of the columns of
$\Phi$. To extract the information about~$x$ that $b$ contains, we use a
\emph{decoder} $\Delta: \R^m \rightarrow \R^n$ which is a (not necessarily
linear) mapping. Then $y \define \Delta(b) = \Delta(\Phi x)$ is our
approximation to $x$ from the information given in $b$. In general, for a
given $b$ and matrix $\Phi$, $\Delta(b)$ may not be unique and it could be
a set of vectors. But here for simplicity we take one of them and deal with
this vector only.

Let $\Sigma_k \define \{z \in \R^n \suchthat \norm{z}{0} \leq k\}$ denote
the vectors in $\R^n$ of support less than or equal to $k$. In the
following we use the classical $\ell_p$-norm, but also other norms are
possible, see Theorem~\ref{thm:CohenBestK1} below. We introduce the distance
\[
\sigma_k(x)_{p} \define \min_{z \in \Sigma_k} \norm{x-z}{p},
\]
and observe that for $x,z \in \R^{n}$ and $p\geq 1$ the following inequality holds:
\begin{equation}\label{eq:SigmaBound}
  \sigma_{2k}(x+z)_p \leq \sigma_k(x)_p + \sigma_k(z)_p.
\end{equation}
We have the following theorem.

\begin{theorem}[\cite{CohDD06a}]\label{thm:CohenBestK1}
  Consider a matrix $\Phi \in \R^{m,n}$ with $m \leq n$, a value $k \in
  \{1, \dots, n\}$, and let $\mathcal{N} = \ker(\Phi)$. If there exists a
  constant $C_0$ such that
  \begin{equation}\label{eq:Sufficond1}
    \norm{\eta}{p} \leq \tfrac{C_0}{2}\, \sigma_{2k}(\eta)_p,\quad \text{for all } \eta \in \mathcal{N},
  \end{equation}
  then there exists a decoder $\Delta$ such that
  \begin{equation}\label{eq:decode}
    \norm{x - \Delta(\Phi x)}{p} \leq C_0\, \sigma_{k}(x)_p, \quad \text{for all } x \in \R^n.
  \end{equation}
  Conversely, if there exists a decoder $\Delta$ such that
  \eqref{eq:decode} holds, then
  \begin{equation}\label{eq:Nesscond1}
    \norm{\eta}{p} \leq C_0\, \sigma_{2k}(\eta)_p, \quad \text{for all } \eta \in \mathcal{N}.
  \end{equation}
\end{theorem}

If we combine Theorem~\ref{thm:CohenBestK1} for $p=1$ with the restricted
isometry property~\eqref{eq:RIPeq1}, then we have the following.

\begin{theorem}[\cite{CohDD06a}]\label{thm:CohenBestK2}
  Let $\Phi \in \R^{m,n}$, $m\leq n$ and $k \in \{1, \dots, n\}$. Assume
  that~$\Phi$ satisfies
  \[
  (1-\delta_{3k}) \norm{x}{2}^2 \leq \norm{\Phi x}{2}^2 \leq
  (1+\delta_{3k}) \norm{x}{2}^2
  \]
  for all $x$ with $\norm{x}{0} \leq 3k$, such that
  \[
  \delta_{3k} \leq \delta < \frac{(\sqrt{2}-1)^2}{3}.
  \]
  Define a decoder $\Delta$ for $\Phi$ via
  \[
  \Delta(b) \define \argmin \{\norm{x}{1} \suchthat b = \Phi x\}.
  \]
  Then
  \[
  \norm{x-\Delta(\Phi x)}{1} \leq C_0\, \sigma_{k}(x)_{1},
  \]
  where
  \[
  C_0 = 2\frac{\sqrt{2} + 1 - (\sqrt{2}-1) \delta}{\sqrt{2} - 1 -
    (\sqrt{2} + 1)\delta}.
  \]
\end{theorem}

Theorem~\ref{thm:CohenBestK2} shows that the $\ell_1$-norm solution can be
as good as the best $k$-term approximation. An analogous result is the
following.

\begin{theorem}[\cite{Can08}]\label{thm:CandesTao}
  Let $\Phi \in \R^{m,n}$, $m\leq n$ and $k \in \{1, \dots, n\}$. Assume
  that $\Phi$ satisfies the restricted isometry property~\eqref{eq:RIPeq1}
  of order $2k$ such that $\delta_{2k}< \sqrt{2}-1$ and $b=\Phi x+e$ where $\norm{e}{2}\leq \epsilon$. 
 If
  \[
  \Delta(b) = \argmin \{\norm{z}{1} \suchthat \| b - \Phi z \|_2 \leq \epsilon\},
  \]
  then
  \[
  \norm{x-\Delta(b)}{2} \leq C_1\, \frac{\sigma_k(x)_{1}}{\sqrt{k}}+C_2 \, \epsilon
  \]
  for some constants $C_1$ and $C_2$ only depending on $\delta_{2k}$.
\end{theorem}
\begin{remark}
It is easy to see that for the case where $\epsilon=0$ and $x$ is $k$-sparse, we have the exact recovery, in other words, $x=\Delta(b)$; see~\cite{Can08} for details.
\end{remark}

Besides the $k$-restricted isometry constant $\delta_k$, a second quantity
plays an important role in compressed sensing~\cite{DonH01,Tro04,Tro06a}.

\begin{definition}
  Let $\Phi \in \R^{m,n}$ with $m \leq n$ have unit norm columns, i.e.,
  $\Phi = [\phi_1 \cdots \phi_n]$ with $\norm{\phi_i}{2} = 1$, for $i = 1,
  \dots, n$. Then the \emph{mutual incoherence} of the matrix $\Phi$ is
  defined by
  \[
  \MutInc{\Phi} := \max_{i \neq j} \abs{\iprod{\phi_i}{\phi_j}}.
  \]
\end{definition}

The mutual incoherence $\MutInc{\Phi}$ of a matrix $\Phi$ is related to the
$k$-restricted isometry constant via
\[
\delta_k \leq (k-1) \MutInc{\Phi}.
\]

The following Lemma shows how the mutual incoherence may be used to bound
the norm of the encoded vector $b = \Phi x$.

\begin{lemma}\label{lemma:Mutual1}
  Let $\Phi = [\phi_1 \cdots \phi_n] \in \R^{m,n}$ with $m \leq n$ have
  unit norm columns. Then for every $x \in \R^n$ the inequality
  \[
  \norm{\Phi x}{2}^2 \leq \big(1 - \MutInc{\Phi}\big) \norm{x}{2}^2 + \MutInc{\Phi}
  \norm{x}{1}^2.
  \]
  holds.
\end{lemma}

\begin{proof}
  The proof follows by the following (in)equalities.
  \begin{align*}
    \norm{\Phi x}{2}^2 & = \sum_{i=1}^n \sum_{j=1}^n x_i \, x_j \iprod{\phi_i}{\phi_j}
    = \norm{x}{2}^2 + \sum_{i \neq j} x_i \, x_j \iprod{\phi_i}{\phi_j} \\
    & \leq \norm{x}{2}^2 + \MutInc{\Phi} \sum_{i \neq j} \abs{x_i}\, \abs{x_j}
    = \norm{x}{2}^2 + \MutInc{\Phi} \big(\norm{x}{1}^2 - \norm{x}{2}^2\big).
  \end{align*}
\end{proof}

Lemma~\ref{lemma:Mutual1} states that $\norm{\Phi x}{2}^2$ is bounded by a
convex combination of $\norm{x}{2}^2$ and $\norm{x}{1}^2$ with the mutual
incoherence as a parameter.

Compressed Sensing (Compressive Sampling) refers to a problem of
``efficient'' recovery of an unknown vector $x\in \R^n$ from the
\emph{partial} information provided by linear measurements $\langle x,
\phi_j \rangle, \phi_j\in\R^n, j = 1,\ldots , m$.  The goal in compressed
sensing is to design an algorithm that approximates $x$ from the
information $b=(\langle x, \phi_1 \rangle, \ldots, \langle x, \phi_m
\rangle) \in \R^m$.  Clearly the most important case is when the number of
measurements $m$ is much smaller than $n$. The crucial step for this to
work, is to build a set of sensing vectors $\phi_j\in \R^n,\ j = 1,\ldots,
m$ that is ``good'' for the approximation of all vectors $x\in \R^n$.
Clearly, the terms ``efficient'' and ``good'' should be clarified in a
mathematical setting of the
problem. 

A natural variant of this setting, and this is the approach that is
discussed here, uses the concept of sparsity. The problem can then be
stated as follows.  For given integers $m\leq n$ we want to determine the
largest sparsity $k(m, n)$ such that there exists a set of vectors
$\phi_j\in \R^n, j = 1,\ldots, m$ and an efficient decoder~$\Delta$,
mapping $b$ into $\R^n$ in such a way that for any $x$ of sparsity $k(m,
n)$ one has exact recovery $\Delta(x) = x$, (see~\cite{Don06a}).

\section{Sparse Representations of Solutions of PDEs}
\label{sec:PDEApproach}

As discussed in the introduction, we want to use similar ideas as those used
in compressed sensing in the context of the solution of partial differential
equations.

\subsection{General Setup}

For a Hilbert space of functions $\HS = \HS(\Omega)$ on a domain $\Omega
\subset \R^d$ with smooth boundary $\Gamma$, we denote by $\hprod{f}{g}$
the inner product of $f, g \in \HS$ and by $\norm{f}{\HS} \define
\sqrt{\hprod{f}{f}}$ the induced norm.

A \emph{generating system} or \emph{dictionary} for~$\HS$ is a family
$\{\phi_{i}\}_{i=1}^{\infty}$ of unit norm elements (i.e.,
$\norm{\phi_{i}}{\HS} = 1$) in $\HS$, such that finite linear combinations
of the elements~$\{\phi_i\}$ are dense in $\HS$. A smallest possible
dictionary is a \emph{basis} of~$\HS$, while the other dictionaries are
\emph{redundant} families of elements. Elements of $\HS$ do not have unique
representations as a linear sums of redundant dictionary elements.

We will consider \emph{elliptic boundary value problems}, see, e.g.,
\cite{Bra07}, and we want to find the solution of
\[
\Lop u = f \quad \text{ in } \Omega,
\]
for a differential operator~$\Lop$ and homogeneous Dirichlet boundary
conditions $u = 0$ on the boundary $\Gamma$ of $\Omega$.

To fix notation and explain our approach, we will lay out a finite element
approach, where we assume that test and solution space $\U \subset \HS$ are
the same. In order to get a closer analogy between the linear algebra
formulation and the function space formulation, we assume that we have a
redundant dictionary $\mathcal{D} = \{\phi_i\}_{i=1}^{\infty}$, such that
\[
\Span \{\phi_j\}_{j=1}^{\infty} = \U.
\]
The corresponding weak formulation for the PDE problem is to find $u \in
\U$ such that
\begin{equation}\label{eq:VarFor}
  a(u,v) \define \hprod{\Lop u}{v} = \hprod{f}{v},\quad\text{for all } v \in \U,
\end{equation}
where $a(\cdot,\cdot)$ is a bilinear form.

Finitely expressing $u,v$ in terms of the dictionary as
\[
u = \sum_{i=1}^{\infty} u_i\, \phi_i \quad\text{and}\quad
v = \sum_{i=1}^{\infty} v_i\, \phi_i,
\]
we can write the problem in terms of infinite vectors and matrices as
\begin{equation}\label{eq:InfQuad}
  \sum_{i,j=1}^{\infty} v_j\, a_{ij}\, u_i = (v^{\infty})\T A^{\infty}
  u^{\infty} = (v^\infty)\T b^{\infty},
\end{equation}
Here, $A^\infty \define [a_{ij}^\infty]$ is the \emph{stiffness matrix},
with $a_{ij}^\infty \define \hprod{\Lop \phi_i}{\phi_j}$, $i,j \in \N$, the
right hand side $b^\infty$ is defined by $b_{i}^\infty \define
\hprod{f}{\phi_i}$, for $i \in \N$, and the coordinate vectors are
\begin{equation}\label{eq:DefUV}
  u^{\infty} \define
  \begin{bmatrix}
    u_1 \\ u_2 \\ \vdots
  \end{bmatrix},\quad
  v^{\infty} \define
  \begin{bmatrix}
    v_1 \\ v_2 \\ \vdots
  \end{bmatrix}.
\end{equation}
The weak solution then can be formulated as the (infinite) linear system
\begin{equation}\label{eq:InfLin}
  A^{\infty} u^{\infty} = b^{\infty}.
\end{equation}
Note that if the $\phi_j$'s are not linearly independent, the (infinite)
linear system~\eqref{eq:InfQuad} and hence~\eqref{eq:InfLin} may not be
solvable or may have many solutions.

\subsection{Algorithmic Framework}
\label{sec:Algorithm}

Let us now consider finite dimensional subproblems of (\ref{eq:InfLin}) by
assuming that the dictionary $\mathcal{D}$ subsumes a \emph{refinement}
procedure, i.e., for each basis function~$\phi_j$ there exists a set of
refined basis functions in $\mathcal{D}$. In particular, we assume that
there is a mapping of each basis function $\phi_i$ to the index set
$\Ref{i} \subset \N$ of refined basis functions. In a \emph{hierarchical}
refinement, every~$\phi_i$ can be written as a linear combination of
$\{\phi_j\}_{j \in \Ref{i}}$.

In many practical applications, the refinement will arise from a geometric
refinement, for instance, by subdividing some triangulation and
corresponding basis functions. Furthermore, the refinement may satisfy $\Ref{j} \cap \Ref{\ell} = \varnothing$ for $j \neq \ell$ (but see Section~\ref{sec:numerics}). For
notational convenience, we define
\[
\Ref{S} \define \bigcup_{j \in S} \Ref{j},
\]
for $S \subset \N$. With this notation, we obtain a sequence of index sets
\[
T^0 \define \{1\},\; T^1 = \Ref{T^0},\; T^2 = \Ref{T^1}, \dots
\]
Define $S^k \define T^0 \cup \dots \cup T^k$ and denote the corresponding
nested subspaces as
\[
\U^0 \subset \U^1 \subset \U^2 \subset \dots \subset \U
\qquad\text{with}\qquad
\U^k \define \Span \{ \phi_j \}_{j \in S^k}.
\]

We will appropriately select subsets $R^k \subseteq C^k \subseteq S^k$,
where $R^k$ and~$C^k$ are seen as subsets of the rows and columns of
$A^\infty$, respectively. The corresponding submatrix is defined as
follows.
\[
A^k \define A^\infty[R^k,C^k] \define [\hprod{\Lop \phi_i}{\phi_j}]_{i \in
  R^k, j \in C^k}.
\]
The corresponding right hand side is
\[
b^k \define b^{\infty}[R^k] \define [\hprod{f}{\phi_j}]_{j \in R^k}.
\]
We thus arrive at the finite dimensional subsystem
\begin{equation}\label{eq:IterISystem}
  A^k x^k = b^k,
\end{equation}
and the approximate solution in this case is
\[
u^k \approx \sum_{j \in C^k} x_j^k\, \phi_j.
\]
As in classical adaptive methods the hope is to keep the size of the matrix
$A^k$ small (i.e., keep $R^k$ and $C^k$ small) and still obtain a good
approximation of the solution arising from the full refinement of
level~$k$, i.e, a solution obtained from solving~\eqref{eq:IterISystem} for
$C^k = R^k = S^k$.

We will now explain how compressed sensing can be used to select
small~$R^k$ and $C^k$ under the condition that we still obtain a convergent
method. For this, assume that we have already selected $R^{k-1}$ and $C^{k-1}$. We
may start with $R^0 = C^0 = \{1\}$, but in practice one should choose an
appropriately fine level. We now refine these sets to $\hat{R}^k \define
\Ref{R^{k-1}}$ and $\hat{C}^k \define \Ref{C^{k-1}}$. Then~$A^k$ and $b^k$ can be
partitioned as follows
\begin{equation}\label{eq:Partition}
  A^k =
  \begin{bmatrix}
    A_{11} & A_{12} & A_{13} \\
    A_{21} & A_{22} & A_{23} \\
    A_{31} & A_{32} & A_{33}
  \end{bmatrix},
  \qquad
  b^k =
  \begin{bmatrix}
    b_{1} \\  b_2 \\ b_3
  \end{bmatrix},
\end{equation}
where
\begin{align*}
  A_{11} & = A^\infty[R^{k-1},C^{k-1}] &
  A_{12} & = A^\infty[R^{k-1},\hat{C}^k] &
  A_{13} & = A^\infty[R^{k-1},\overline{C}^{k}]\\
  A_{21} & = A^\infty[\hat{R}^k,C^{k-1}] &
  A_{22} & = A^\infty[\hat{R}^k,\hat{C}^k] &
  A_{23} & = A^\infty[\hat{R}^k,\overline{C}^{k}]\\
  A_{31} & = A^\infty[\overline{R}^{k},C^{k-1}] &
  A_{32} & = A^\infty[\overline{R}^{k},\hat{C}^k] &
  A_{33} & = A^\infty[\overline{R}^{k},\overline{C}^{k}],
\end{align*}
with $\overline{R}^k \define \N \setminus (R^{k-1} \cup \hat{R}^k)$ and
$\overline{C}^k$ defined analogously. Similarly, the right hand side is
defined as
\[
b_1 = b^\infty[R^{k-1}],\quad
b_2 = b^\infty[\hat{R}^k],\quad
b_3 = b^\infty[\overline{R}^k].
\]

\begin{figure}
  \includegraphics[width=0.45\textwidth]{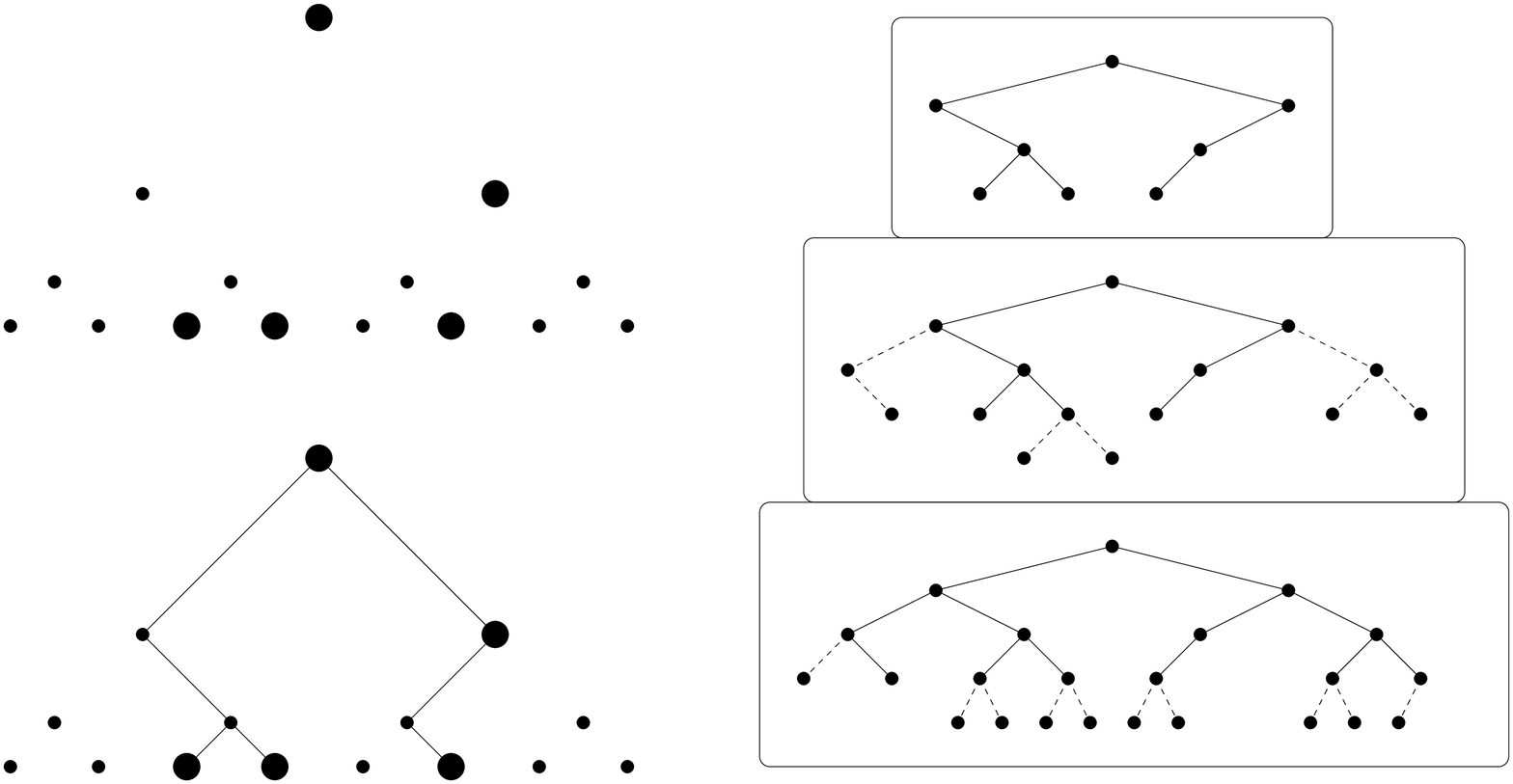}
  \hfill
  \psfrag{I}{}
  \includegraphics[width=0.45\textwidth]{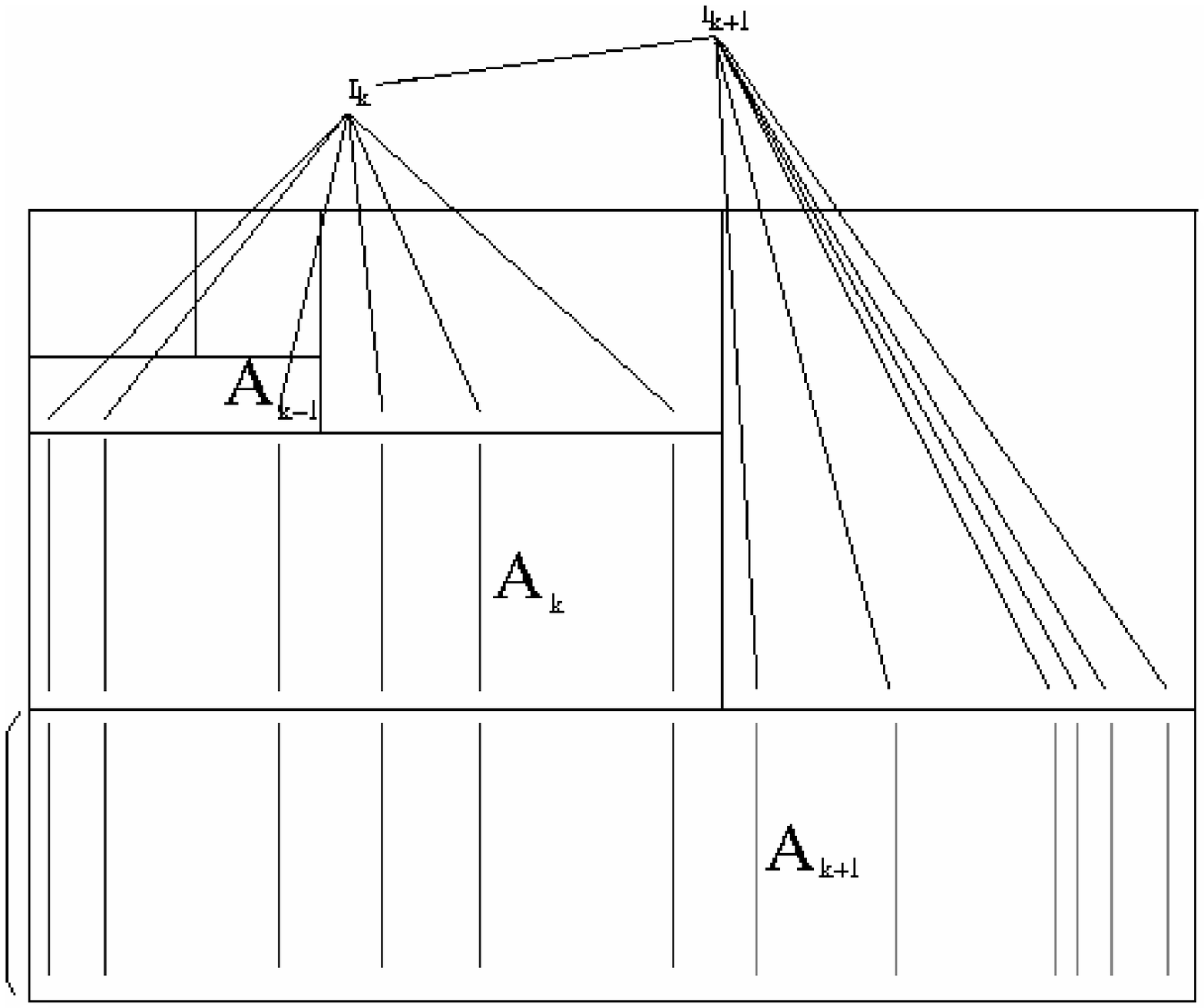}
  \caption{Illustration of the stepwise refinement method. \emph{Left:}
    Picture of a tree that is iteratively refined. \emph{Right:} A view on
    the matrix $A^\infty$ and its partition as in~\eqref{eq:Partition}.}
  \label{fig:RefinementIllustration}
\end{figure}

The main idea in the construction of $R^{k}$ and $C^{k}$ is to use the
underdetermined subsystem
\[
[A_{21}\; A_{22}]\, z = b_2,
\]
and to compute a sparse solution, by taking a minimal $\ell_1$-solution,
i.e.,
\[
z^k = \argmin \{ \norm{z}{1} \suchthat [A_{21}\; A_{22}]\, z = b_{2}\}.
\]
For the noisy case, we solve the following problem:
\[
z^k = \argmin \{ \norm{z}{1} \suchthat \norm{ [A_{21}\; A_{22}]\, z - b_{2} }{2} \leq \epsilon_k\}.
\]
where $\epsilon_k$ is a given upper bound on the size of the noise.
The submatrix $[A_{21}\; A_{22}]$ is chosen because it combines the
refined rows with the full set of columns that are available at the current
iteration.

Now assume that $C \subseteq C^{k-1} \cup \hat{C}^k$ is the index set
corresponding to the support of~$z^k$ as defined
in Section~\ref{sec:not}. Then the new sets are set to
\[
C^{k} = C \cup C^{k-1},\qquad R^{k} = C^{k}.
\]
Thus, the support of $z^k$ is only used to select basis functions among
$\hat{C}^k$ and the information in $C^k \cap C$ is not used.

Note that by construction
\[
C^{k} \subseteq \Tree{C},
\]
where $\Tree{C}$ is the set of basis functions on the path of a basis
function to the root of the refinement-tree, i.e.,
\[
\Tree{j} \define \{ \ell \in \N \suchthat \exists j_1, \dots, j_s \text{
  with } j_1 \in \Ref{\ell},\; j_2 \in \Ref{j_1},\; \dots,\; j \in
\Ref{j_s}\}.
\]
The process is terminated if $\norm{z^{k} - z^{k-1}}{2} \leq \varepsilon$,
where $\varepsilon$ is a given tolerance.  Since we are using
$\ell_1$-mini\-miz\-ation, which in~\cite{Che95} is called basis pursuit,
we call this process \emph{Iteratively Refined Basis Pursuit.}  The method
is summarized as Algorithm~\ref{alg:IRBP}.
Figure~\ref{fig:RefinementIllustration} gives an illustration of the
process.

\begin{algorithm}[tb]
  \caption{Iteratively Refinement Basis Pursuit (IRBP)}
  \label{alg:IRBP}
  \begin{algorithmic}[1]
    \STATE Set $R^0 = C^0 = \{1\}$
    \FOR{$k = 1, \dots, $ until convergence}
    \STATE Construct $\hat{R}^k = \Ref{R^{k-1}}$, $\hat{C}^k = \Ref{C^{k-1}}$
    \STATE Construct $A_{21} = A^\infty[\hat{R}^k,C^{k-1}]$, $A_{22} =
    A^\infty[\hat{R}^k,\hat{C}^k]$, and $b_2 = b^\infty[\hat{R}^k]$.
    \STATE Solve the following minimization problem:
    \[
    z^k = \argmin \{\norm{z}{1} \suchthat [A_{21}\; A_{22}] z = b_2\}.
    \]
    \STATE Let $C \subseteq C^{k-1} \cup \hat{C}^k$ be the index set
    corresponding to the support of $z^k$.
    \STATE Set $C^{k} = C \cup C^{k-1}$, $R^{k} \define C^{k}$.
    \ENDFOR
  \end{algorithmic}
\end{algorithm}

\begin{figure}
  \centering
  \psfrag{1}{1}
  \psfrag{2}{2}
  \psfrag{3}{3}
  \psfrag{4}{4}
  \psfrag{5}{5}
  \psfrag{6}{6}
  \psfrag{7}{7}
  \includegraphics[width=0.4\textwidth]{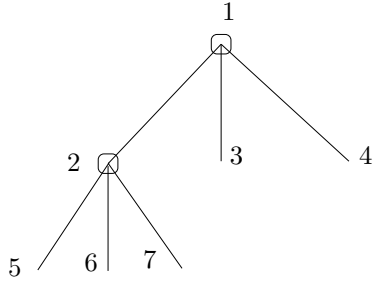}
  \caption{Illustration of Example~\ref{ex:Refinement}.}
  \label{fig:TreeExample}
\end{figure}

\begin{example}\label{ex:Refinement}
  By definition, the first sets are $R^0 = C^0 = \{1\}$.  Now suppose that
  $\Ref{1} = \{2,3,4\}$, i.e, the initial rows and columns are $\hat{C}^1 =
  \hat{R}^1 = \{2,3,4\}$. We then solve the $\ell_1$-mini\-miz\-ation
  problem for the corresponding $3 \times 4$ matrix. Now suppose that
  $\supp{z^1} = \{1,2\}$, then $C^1 = \{1,2\}$. Assume that $\Ref{2} =
  \{5,6,7\}$, see Figure~\ref{fig:TreeExample} for an illustration. Then
  $\hat{C}^2 = \{3,4,5,6,7\}$ and the next matrix is of size $5 \times
  7$, since $R^1 = \{1,2\}$ and $\hat{R}^2 = \{3,4,5,6,7\}$.
\end{example}

Note that if the support of $z^k$ is full, then the above procedure yields a
simple refinement process in which $R^k = S^{k-1}$ and $C^k = S^k$.

In general, this process need not converge. In fact, if at some level we
see that the error does not decrease, then we back up in the tree of
refined basis functions and refine at higher levels until we obtain a
decrease in the error.

\begin{remark}
  Note that although our description was based on the assumption that test
  and solution space are the same, the principle of the process does not
  depend on this assumption. Similar concepts for adaptive refinement
  methods in the context of wavelets are presented, e.g.
  in~\cite{Bar99,CohDDD01,CohDD01a}.
\end{remark}

\subsection{Properties of the Proposed Method}

In our approach we want to achieve several goals. The solution~$z^k$ should
be sparse, and~$z^k$ should be a good approximation of the solution
$x^{i+1}$ of~\eqref{eq:IterISystem}. In order to analyze the behavior of
the suggested approach, we study the case of two levels and assume that
$R^k = C^k$. We will also slightly change notation as follows. For $k \in
\N$, we set $n \define \card{R^{k-1}}$ and $N \define \card{(R^{k-1} \cup
  \Ref{R^{k-1}})}$. We then introduce the following notation for
submatrices of $A^\infty$ and subvectors of $b^\infty$ as
in~\eqref{eq:Partition}:
\begin{equation}\label{eq:DefAN}
  A^N =
  \begin{bmatrix}
    A_{11} & A_{12} \\
    A_{21} & A_{22}
  \end{bmatrix},
  \qquad
  b^N =
  \begin{bmatrix}
    b_{1} \\  b_2
  \end{bmatrix},
\end{equation}
and $A^n \define A_{11}$, $b^n \define b_1$. Note that $A^N$ is of size $N
\times N$ and $A^n$ of size $n \times n$. The corresponding linear systems
are
\begin{equation}\label{eq:InfLin2}
  A^{\infty} x^{\infty} = b^{\infty},
\end{equation}
\begin{equation}\label{eq:FinLinN}
  A^N x^N = b^N,
\end{equation}
and $A^n x^n = b^n$. In Algorithm~\ref{alg:IRBP}, we take the last $N-n$
rows of~\eqref{eq:DefAN} and consider the linear system
\begin{equation}\label{eq:AN-n}
  A^{N-n,N} z^N = b_{2}\qquad\text{ with }\qquad
  A^{N-n,N} \define
  [A_{21}\, A_{22}].
\end{equation}
We find a solution of this underdetermined system, by taking a minimal
$\ell_1$-solution, i.e., by solving
\begin{equation}\label{eq:Minl1}
  z^N = \argmin \{ \norm{z}{1} \suchthat A^{N-n,N} z = b_{2}\}.
\end{equation}

As a first step of the analysis we estimate the energy norm error
between~$x^N$ and $z^N$ in terms of the difference
$\norm{x^N-x^{\infty}}{}$. To derive such a bound, we need to embed $z^N$
and $x^n$ into $\R^\infty$ by appending $0$ as follows:
\begin{equation}\label{eq:HatNot}
  \hat{z}^N = \begin{bmatrix} z^N \\ 0 \end{bmatrix},\
  \hat{x}^N = \begin{bmatrix} x^N \\ 0 \end{bmatrix}.
\end{equation}
We assume that $A^{\infty} x^{\infty} = b^{\infty}$, $A^N x^N = b^N$, and
that $z^N$ is determined by~\eqref{eq:Minl1}.

We want to find necessary and sufficient conditions on the
matrix~$A^{\infty}$ and thus on the dictionary $\{\phi_i\}$, such that
there exists a constant $C_{\frac{n}{N}}$ for which the following
inequality holds,
\begin{equation}\label{eq:ENest}
  (\hat{z}^N - \hat{x}^N)\T A^{\infty} (\hat{z}^N - {\hat{x}}^N)
  \leq C_{\frac{n}{N}}\,
  (x^{\infty} - \hat{x}^N)\T A^{\infty} (x^{\infty}-\hat{x}^N).
\end{equation}
The constant $C_{\frac{n}{N}}$ should only depend on the ratio
$\frac{n}{N}$. Considering the results of~\cite{CanT06,CohDD06a,DonH01}, we
may expect that if the matrices $A^\infty$, $A^N$ have a small mutual
incoherence or some restricted isometry property, such conditions can be
obtained. We will come back to this point in Section~\ref{sec:rippde}.

If Inequality~\eqref{eq:ENest} holds, then by using the triangle inequality
we obtain the error estimate
\[
  (x^{\infty} - \hat{z}^N)\T A^{\infty} (x^{\infty} - \hat{z}^N)
  \leq (1 + C_{\frac{n}{N}})
  (x^{\infty} - \hat{x}^N)\T A^{\infty} (x^{\infty}-\hat{x}_N),
\]
which means that the (hopefully sparse) solution $z^N$ obtained by solving
\eqref{eq:Minl1} is as almost good as the solution of~\eqref{eq:FinLinN}.

To estimate the errors between $x^\infty$, $x^N$, and $z^N$, we need to
consider a refined partition of $A^\infty$ and $b^\infty$ as defined
in~\eqref{eq:DefAN}. In order to relate the solutions at different levels
of refinement and the solution of the $\ell_1$-minimization, we make use of
the following Lemmas.

\begin{lemma}\label{Lem:ErrorEst1}
  Let
  \[
  x^{\infty} =
  \begin{bmatrix}
    x^{\infty}_{1} \\
    x^{\infty}_{2} \\
    x^{\infty}_{3}
  \end{bmatrix},\
  x^N =
  \begin{bmatrix}
    x^{N}_{1} \\
    x^{N}_{2}
  \end{bmatrix}, \text{and }
  z^N =
  \begin{bmatrix}
    z^{N}_{1} \\
    z^{N}_{2}
  \end{bmatrix}
  \]
  be solutions of the problems~\eqref{eq:InfLin2}, \eqref{eq:FinLinN},
  and~\eqref{eq:Minl1}, respectively. Furthermore, let~$\hat{x}^N$ and
  $\hat{z}^N$ be as in~\eqref{eq:HatNot}. Then inequality~\eqref{eq:ENest}
  can be rewritten as
  \[
  (x^{N}_{1} - z^{N}_{1})\T \big(b_1 - [A_{11}\, A_{12}]
  \begin{bmatrix}
    z^{N}_{1} \\
    z^{N}_{2}
  \end{bmatrix}
  \big)
  \leq  C_{\frac{n}{N}}\; (x^{\infty}_{3})\T
  \big(b_3 - [A_{31}\, A_{32}]
  \begin{bmatrix}
    x^{N}_{1} \\
    x^{N}_{2}
  \end{bmatrix}
  \big).
  \]
\end{lemma}
\begin{proof}
  Since
  \[
  A^{\infty} x^{\infty} = b^{\infty},\ A^{\infty} \hat{x}^N =
  \begin{bmatrix}
    b_1 \\
    b_2 \\
    [A_{31}\, A_{32}] x^N
  \end{bmatrix},\
  A^{\infty} \hat{z}^N =
  \begin{bmatrix}
    [A_{11}\, A_{12}] z^N \\
    b_2 \\
    [A_{31}\, A_{32}] z^N
  \end{bmatrix},
  \]
  it follows that
  \[
  A^{\infty} (\hat{x}^N - x^{\infty}) =
  \begin{bmatrix}
    0 \\
    0 \\
    [A_{31}\, A_{32}] z^N - b_3
  \end{bmatrix}
  \]
  and
  \[
  A^{\infty} (\hat{z}^N - \hat{x}^N) =
  \begin{bmatrix}
    [A_{11}\, A_{12}] z^N - b_1 \\
    0 \\
    [A_{31}\, A_{32}] (z^N - x^N)
  \end{bmatrix}.
  \]
  Thus, we have
  \[
  (\hat{x}^N - x^{\infty})\T A^{\infty} (\hat{x}^N - x^{\infty}) =
  (x_3^{\infty})\T \big(b_3 - [A_{31}\, A_{32}]
  \begin{bmatrix}
    x^N_1 \\
    x^N_2
  \end{bmatrix}
  \big)
  \]
  and
  \[
  (\hat{z}^N - \hat{x}^N)\T A^{\infty} (\hat{z}^N - \hat{x}^N) =
  (x^{N}_{1} - z^{N}_{1})\T \big(b_1 - [A_{11}\, A_{12}]
  \begin{bmatrix}
    z^{N}_{1} \\
    z^{N}_{2}
  \end{bmatrix}
  \big).
  \]
  Plugging these expressions into~\eqref{eq:ENest} yields the claim.
\end{proof}

\begin{remark}\label{RemConj1}
  A weaker version of (\ref{eq:ENest}) and of Lemma \ref{Lem:ErrorEst1}
  will be given in Section~\ref{sec:rippde}.
\end{remark}

The following Lemma gives a condition that has to be satisfied in order to
guarantee that the refinement process can be iterated.

\begin{lemma}\label{lemma:Norml1}
  Let
  \[
  z^N =
  \begin{bmatrix}
    z^{N}_{1} \\ z^{N}_{2}
  \end{bmatrix}
  \]
  be a solution of~\eqref{eq:Minl1}, where $A^{N-n,N}$
  is as defined in~\eqref{eq:AN-n}, and suppose that $A_{22}$ is
  invertible.  If $z^N_1 \neq 0$, then
  \begin{equation}\label{eq:NormIneql1}
    \norm{A^{-1}_{22} A_{21}}{1} \geq 1.
  \end{equation}
\end{lemma}

\begin{proof}
  Since
  \[
  z' =
  \begin{bmatrix}
    0 \\
    A^{-1}_{22} b_2
  \end{bmatrix}
  \]
  is a feasible solution of~\eqref{eq:Minl1}, it follows that
  \begin{equation}\label{eq:MinProp}
    \norm{z^N}{1} \leq \norm{z'}{1}.
  \end{equation}
  Moreover, from $A_{21} z^{N}_{1} + A_{22} z^{N}_{2} = b_2$, we obtain
  that $z^{N}_{2} = A^{-1}_{22} b_2 - A^{-1}_{22} A_{21} z^{N}_{1}$. Thus,
  using~\eqref{eq:MinProp}, we obtain
  \begin{align*}
    \norm{A^{-1}_{22}b_2}{1} & \geq \norm{z^{N}}{1}\\
    & = \norm{z^{N}_{2}}{1} + \norm{z^{N}_{1}}{1} \\
    & = \norm{A^{-1}_{22} b_2 - A^{-1}_{22} A_{21} z^{N}_{1}}{1} + \norm{z^{N}_{1}}{1}\\
    & \geq \norm{A^{-1}_{22}b_2}{1} - \norm{A^{-1}_{22} A_{21}
      z^{N}_{1}}{1} + \norm{z^{N}_{1}}{1}.
  \end{align*}
  Since $\norm{A^{-1}_{22} A_{21} z^{N}_{1}}{1} \leq \norm{A^{-1}_{22}
    A_{21}}{1}\cdot \norm{z^{N}_{1}}{1}$, we have
  \[
  \norm{z^{N}_{1}}{1} \leq \norm{A^{-1}_{22} A_{21}}{1}\cdot \norm{z^{N}_{1}}{1},
  \]
  which completes the proof.
\end{proof}

\begin{remark}
  Lemma~\ref{lemma:Norml1} implies that a solution of the
  $\ell_1$-minimization problem can only lead to an improvement if the
  matrix $[A_{21}\, A_{22}]$ satisfies~\eqref{eq:NormIneql1}. Otherwise, an
  optimal solution can already be obtained by solving the linear system
  $A_{22} z^N_2 = b_2$. Another observation is that
  Lemma~\ref{lemma:Norml1} remains true for any nonsingular principal
  submatrix of $A^N$.
\end{remark}

In order to compare sparse and non-sparse solutions we introduce the short
notation $s(x) \define \supp{x}$ for a vector $x \in \R^m$ and
$\overline{s}(x) = \{1,2,\dots,m\} \setminus s(x)$. For $x \in \R^m$ and $S
\subset \{1, \dots, m\}$, we denote
\[
(y_S)_i =
\begin{cases}
  y_i & \text{if }i \in S\\
  0   & \text{otherwise}.
\end{cases}
\qquad
i = 1, \dots, m.
\]

We then have the following Lemma.

\begin{lemma}\label{lemma:Norml1-2}
  Let $z^N$ be a solution of~\eqref{eq:Minl1}, where $A^{N-n,N}$ is as
  in~\eqref{eq:AN-n}, and let $x^N$ be a solution of $A^N x^N = b^N$, with
  $A^N$ as in~\eqref{eq:DefAN}. Then for the difference $\delta^N \define
  z^N - x^N$ we have the inequality
  \[
  \frac{\norm{\delta^N_{s(x^N)}}{1}}{\norm{\delta^N}{1}} \geq \frac{1}{2}.
  \]
\end{lemma}

\begin{proof}
  Since $x^N$ is a feasible solution of~(\ref{eq:Minl1}), we have
  \begin{equation}\label{eq:MinProp2}
    \norm{x^N + \delta^N}{1} = \norm{z^N}{1} \leq \norm{x^N}{1}.
  \end{equation}
  Furthermore, $x^N + \delta^N = x^N + \delta^N_{s(x^N)} +
  \delta^N_{\overline{s}(x^N)}$. Therefore, by~\eqref{eq:MinProp2}, we have
  \begin{align*}
    \norm{x^N}{1} & \geq \norm{x^N+\delta^N}{1} \\
    & = \norm{x^N + \delta^N_{s(x_N)}}{1} + \norm{\delta^N_{\overline{s}(x^N)}}{1} \\
    & \geq \norm{x^N}{1} - \norm{\delta^N_{s(x^N)}}{1} + \norm{\delta^N_{\overline{s}(x^N)}}{1}.
  \end{align*}
  Rewriting yields that
  \[
  \norm{\delta^N_{s(x^N)}}{1} \geq \norm{\delta^N_{\overline{s}(x^N)}}{1} =
  \norm{\delta^N - \delta^N_{s(x^N)}}{1}
  \geq \norm{\delta^N}{1} - \norm{\delta^N_{s(x^N)}}{1},
  \]
  which implies the assertion.
\end{proof}

\begin{remark}
  The proof of Lemma~\ref{lemma:Norml1-2} shows that
  \[
  \norm{(z^N - x^N)_{\overline{s}(x^N)}}{1} \leq \norm{(z^N - x^N)_{s(x^N)}}{1}.
  \]
  In particular, if $(z^N - x^N)_{s(x^N)} = 0$, then we conclude that $z^N
  = x^N$.  If instead of $\ell_1$-minimization, we use
  $\ell_0$-minimization and compute
  \[
  w^N = \argmin \{ \norm{z}{0} \suchthat A^{N-n,N} z = b_2\},
  \]
  then we get the analogous estimate
  \[
  \norm{(w^N-x^N)_{\overline{s}(x^N)}}{0} \leq \norm{(w^N -
    x^N)_{s(x^N)}}{0}.
  \]
\end{remark}

\begin{remark}
  In general, it is not true that $z^N$ and $x^N$ satisfy the inequality
  $\norm{z_N}{0} \leq \norm{x_N}{0}$.  For example if
  \[
  A^4 = \left[
    \begin{array}{cccc}
      1 & 0 & \frac{1}{\sqrt{2}} & 0   \\
      0 & 1 & -\frac{1}{2} & 0  \\
      \frac{1}{\sqrt{2}} & -\frac{1}{2} & 1 & -\frac{1}{2} \\
      0 & 0 & -\frac{1}{2} & 1   \\
    \end{array}
  \right]
  \]
  then
  \[
  A^{3,4} =
  \left[
    \begin{array}{cccc}
      0 & 1 & -\frac{1}{2} & 0  \\
      \frac{1}{\sqrt{2}} & -\frac{1}{2} & 1 & -\frac{1}{2} \\
      0 & 0 & -\frac{1}{2} & 1   \\
    \end{array}
  \right]
  \]
  Now consider $b^4=[\sqrt{2}, 6, -\frac{3}{2}, -1]\T$ and
  $x^4 = [0, 7, 2, 0]\T$. Then
  \[
  z^4=[\sqrt{2}, 6, 0, -1]\T,
  \]
  and in this case $\norm{z_N}{1} = 7+\sqrt{2} < 7+2 = \norm{x_N}{1}$ and
  \[
  \norm{z_N}{0}=3 > 2=\norm{x_N}{0}.
  \]
\end{remark}

In this section we have set the stage for the solution of PDEs via
$\ell_1$-mini\-miz\-ation. In the next section we provide details.

\section{The Restricted Isometry Property for Elliptic PDEs}
\label{sec:rippde}

In this section, we again discuss the special case that solution space and
test space are the same, i.e., we assume $\U = \V \subset \HS$ and we
consider the symmetric bilinear form $a(\cdot,\cdot): \U \times \V
\rightarrow \R$ associated with the operator~$\Lop$ as
in~\eqref{eq:VarFor}.

We also assume that there exist constants $\alpha_1 > 0$, $\alpha_2 <
\infty$, such that:
\begin{equation}\label{eq:BoundUnifEll}
  \alpha_1\, \norm{u}{\HS}^2 \leq a(u,u) \leq \alpha_2\, \norm{u}{\HS}^2,
\end{equation}
i.e., $a(\cdot,\cdot)$ is \emph{uniformly elliptic} with constant
$\alpha_1$ and \emph{uniformly bounded} with constant $\alpha_2$.

In order to connect this classical norm equivalence with the $k$-restricted
isometry property, we assume that for the dictionary $\mathcal{D} =
\{\phi_k\}_{i=1}^{\infty}$ the following \emph{$k$-equivalence} between
$\norm{\cdot}\HS$ and the $\ell_2$-norm $\norm{\cdot}{2}$ holds, i.e., we
assume that there exist constants $\beta_1 > 0$, $\beta_2 < \infty$ with
\begin{equation}\label{eq:NormEquivalence}
  \beta_1\, \norm{u^\infty}{2} \leq
  \Norm{\sum_{i=1}^\infty u_i\, \phi_i}{\HS} \leq
  \beta_2\, \norm{u^\infty}{2},
\end{equation}
for all infinite vectors $u^\infty$ as in~\eqref{eq:DefUV} with the
property that $\norm{u^\infty}{0} \leq k$.

Note that inequality~\eqref{eq:BoundUnifEll} can be written as
\[
  \alpha_1\, \Norm{\sum_{i=1}^\infty u_i\, \phi_i}{\HS}^2 \leq
  \norm{{{A^{\infty}}^\frac{1}{2}}\, u^\infty}{2}^2 \leq
  \alpha_2\, \Norm{\sum_{i=1}^\infty u_i\, \phi_i}{\HS}^2,
\]
or equivalently as
\begin{equation}\label{eq:RIPcond}
  \alpha_1\, \Norm{\sum_{i=1}^\infty u_i\, \phi_i}{\HS}^2 \leq
  (u^\infty)\T A^\infty u^\infty \leq
  \alpha_2\, \Norm{\sum_{i=1}^\infty u_i\, \phi_i}{\HS}^2.
\end{equation}
Combining inequalities~\eqref{eq:NormEquivalence} and~\eqref{eq:RIPcond},
we obtain
\[
\alpha_1\, \beta_1^2\, \norm{u^\infty}{2}^2 \leq
(u^\infty)\T A^{\infty} u^\infty \leq
\alpha_2\, \beta_2^2\, \norm{u^\infty}{2}^2,
\]
for all vectors $u^\infty$ with $\norm{u^\infty}{0} \leq k$.

We consider the dictionary $\mathcal{D} = \{\phi_1, \phi_2, \dots\}$, and
we choose a set $\mathcal{I}^N = \{q_1, \dots, q_N\} \subset \N$ with
associated elements $\phi_{q_1}, \dots, \phi_{q_N} \in \mathcal{D}$. For
the theoretical analysis we may assume w.l.o.g. that $\mathcal{I}^N = \{ 1,
2, \dots, N\}$. This selection can be obtained via an appropriate
reordering of the elements $\phi_i$ of the dictionary.

The corresponding finite stiffness matrix associated with this subset is
then $A^N = [a_{ij}] \in \R^{N,N}$ with $a_{ij} = a(\phi_{i},\phi_{j})$,
$i,j = 1, \dots, N$. Since we have assumed uniform ellipticity and since
test and solution space are equal, it follows that~$A^N$ is symmetric and
positive semidefinite with $\rank(A^N) = N$; $A^N$ is positive definite if
$\phi_1, \dots, \phi_N$ form a basis.

Since $A^N$ is symmetric positive semidefinite, $A^N$ has a factorization
$A^N = (B^N)\T B^N$, where $B^N \in \R^{n,N}$ has full row rank. Hence,
there exists a permutation matrix~$P$ such that
\[
P\, B^N =
\begin{bmatrix}
  B^n \\
  B^{N-n}
\end{bmatrix},
\]
with $B^n \in \R^{n,n}$ invertible. This yields
\[
P\T A^N P=
\begin{bmatrix}
  A_{11} & A_{12}\\
  A_{21} & A_{22}
\end{bmatrix}
\]
with stiffness matrix $A^n = A_{11}$. Then we have
\[
A^N = (P\, B^N)\T P\, B^N =
\begin{bmatrix}
  B^n \\
  B^{N-n}
\end{bmatrix}\T
\begin{bmatrix}
  B^n \\
  B^{N-n}
\end{bmatrix}
= (B^n)\T B^n +  (B^{N-n})\T B^{N-n}.
\]
Suppose that it is possible to choose the permutation matrix $P$ in such a
way that (measured in spectral norm)
\[
\norm{B^{N-n}}{2} \leq \epsilon
\]
with a small $\epsilon > 0$, i.e.,
\[
\norm{B^{N-n} x}{2}^2 \leq \epsilon\, \norm{x}{2}^2
\]
for all $x \in \R^n$. Suppose further that
\[
(1 - \delta_k)\, \norm{x^N}{2}^2 \leq (x^N)\T A^N x^N \leq
(1 + \delta_k)\, \norm{x^N}{2}^2
\]
or equivalently
\[
(1 - \delta_k)\, \norm{x^N}{2}^2 \leq
\norm{B^N x^N}{2}^2 \leq
(1 + \delta_k)\, \norm{x^N}{2}^2
\]
for all $x^N$ with $\norm{x^N}{0}\leq k$.

To get an error estimate between the solution that is based on
$\ell_1$-mini\-miz\-ation and the best $k$-term approximation, we first
prove the following result.

\begin{theorem}\label{thm:ErrorEstimate}
  Let $\mathcal{A} \in \R^{N,N}$ be symmetric positive semi-definite, and
  consider the solvable linear system $\mathcal{A} x = b$. Let $\mathcal{A}
  = B\T B$ be a full rank factorization, and let $P \in \R^{N,N}$ be a
  permutation matrix such that the following properties hold.
  \begin{enumerate}
  \item\label{property:decomp}%
    $P\, B =
    \begin{bmatrix}
      B_1 \\
      B_2
    \end{bmatrix}$, and $\mathcal{A} = B\T B = B_1\T B_1 + B_2\T B_2$;
  \item\label{property:sparse}%
    for any solution $x$ of $\mathcal{A} x = b$, $B_2\T B_2\, x$ is $k$-sparse,
    i.e., $B_2\T B_2\, x \in \Sigma_k$, where $\Sigma_k = \{z \suchthat
    \norm{z}{0} \leq k\}$;
  \item the $2k$-restricted isometry constant $\delta_{2k}$ for $B_1$ is
    sufficiently small (for example $\delta_{2k} < \sqrt{2}-1$).
  \end{enumerate}
  Then
  \begin{equation}\label{eq:EstRIP}
    (x - \tilde{x})\T \mathcal{A} \, (x - \tilde{x}) \leq C_k\, \sigma_k^2(x)_1,
  \end{equation}
where $ \sigma_k(x)_1 = \min_{z \in \Sigma_k}\norm{z-x}{1} $, $\tilde{x}$ is obtained via the solution of the minimization
  problems
  \begin{equation}\label{Eq:DecBP}
  \tilde{y} = \argmin_y \norm{b - B_1\T y}{1}
  \quad\text{ and }\quad
  \tilde{x} = \argmin_x \{ \norm{x}{1}: B_1 x = \tilde{y}\}
  \end{equation}
  and the constant $C_k$ only depends on $k$, the mutual incoherence
  $\MutInc{B_1}$ and $\norm{B_2}{2}$.
\end{theorem}

\begin{proof}
  By Assumption~\ref{property:decomp}, $\mathcal{A} x = b$ implies that $b
  = B_{1}\T B_{1}\, x + B_{2}\T B_{2}\, x$. By
  Assumption~\ref{property:sparse}, it follows that $e = B_{2}\T B_{2}\, x$
  is $k$-sparse. Then by Theorem~1.3 of~\cite{CanT06} we obtain exact
  recovery, i.e., if $B_1\, x = \tilde{y}$, then
  \[
  \tilde{y} = \argmin_y \norm{b - B_1\T y}{1}.
  \]
  The remainder of the proof is then based on
  Theorems~\ref{thm:CohenBestK1}, \ref{thm:CandesTao}, and
  Lemma~\ref{lemma:Mutual1}. Since
  \[
  \mathcal{A} = B\T B = B_1\T B_1 + B_{2}\T B_{2},
  \]
  it follows that
  \begin{equation}\label{eq:Energy1}
    (x - \tilde{x})\T \mathcal{A} (x - \tilde{x}) = \norm{B_1 (x - \tilde{x})}{2}^2 +
    \norm{B_{2} (x - \tilde{x})}{2}^2,
  \end{equation}
  where $\tilde{x} = \argmin_x \{ \norm{x}{1}: B_1 x = \tilde{y}\}$.
By Theorem~\ref{thm:CandesTao}, we also have
  \[
\norm{B_{2}(x - \tilde{x})}{2}^2 \leq \norm{B_2}{2}^2\norm{x - \tilde{x}}{2}^2 \leq  \frac{\norm{B_2}{2}^2 C_{1,k}^2}{k} \, \sigma^2_k(x)_{1},
  \]
  where $C_{1,k}$ only depend on $\delta_k$. 

W.l.o.g. we can assume that $B_1$ has unit norm columns. Otherwise
  instead of the linear equation $B_1x=\tilde{y}$, we can consider the
  following linear equation:
  \[
  (B_1S)S^{-1}x=\tilde{y},
  \]
  where $S=\mathrm{diag}(\frac{1}{\norm{Be_i}{2}})$ and $e_i$ is the $i$-th
  column of the identity matrix. Then the matrix $B_1S$ has unit norm
  column and therefore, by Lemma~\ref{lemma:Mutual1} we have that
  \[
  \norm{B_1(x - \tilde{x})}{2}^2 \leq (1 - \MutInc{B_1})\,
  \norm{x - \tilde{x}}{2}^2 + \MutInc{B_1}\, \norm{x - \tilde{x}}{1}^2.
  \]
  By Theorem~\ref{thm:CandesTao}, we have $\norm{x - \tilde{x}}{2}^2 \leq
  \frac{C_{1,k}^2}{k}\, \sigma^2_k(x)_{1}$ and by
  Theorem~\ref{thm:CohenBestK1} we have $\norm{x - \tilde{x}}{1}^2 \leq
  C_{2,k}^2\, \sigma^2_k(x)_{1}$, where $C_{1,k}$ and $C_{2,k}$ only depend
  on $\delta_k$. Combining these inequalities with~\eqref{eq:Energy1}, we
  get
  \[
  (x - \tilde{x})\T \mathcal{A} (x - \tilde{x}) \leq
  \Big( (1 - \MutInc{B_1}) \frac{C_{1,k}^2}{k} + \MutInc{B_1}\, C_{2,k}^2
  + \frac{\norm{B_2}{2}^2 C_{1,k}^2}{k}\Big) \sigma^2_k(x)_{1}.
  \]
  This concludes the proof.
\end{proof}


Applying Theorem~\ref{thm:ErrorEstimate} to the matrix $\mathcal{A}
=(B^{N-n,N})\T B^{N-n,N}$, where matrix $B^{N-n,N} = [A_{21} A_{22}]$ as
in~\eqref{eq:DefAN}, we obtain the corresponding estimate for the stiffness
matrix.

\begin{corollary}\label{cor:ErrorEstimate}
  Let $\mathcal{A}^N\in \R^{N,N}$ be a symmetric positive semidefinite
  matrix of rank~$N-n$ and $\mathcal{A}^Nx^N=b^N$. Let $\mathcal{A}^N =
  (B^{N-n,N})\T B^{N-n,N}$ be a full rank factorization of $\mathcal{A}^N$,
  where $B^{N-n,N}=[A_{21} \quad A_{22}]$. If the $2k$ restricted isometry
  constant $\delta_{2k}$ for $B^{N-n,N}$ is sufficiently small (e.g. if
  $\delta_{2k}<\sqrt{2}-1$), then
  \[
  (x^N-\hat{x})\T \mathcal{A}^N(x^N-\hat{x})\leq C_k\sigma_k^2(x^N)_1,
  \]
  where $\hat{x}$ is obtained via the solution of the minimization problem
  \[
    \hat{y}=\mathrm{argmin}_{y} \norm{b^N-{B^{N-n,N}}\T y}{1}
  \]
  and
  \[
    \hat{x}=\mathrm{argmin}_{z}\{\norm{z}{1}: B^{N-n,N}z=\hat{y}\},
  \]
  and $C_k$ only depends on $k$ and $\mathcal{M}(B^{N-n,N})$.
\end{corollary}

\begin{proof}
  Taking $B_1=B^{N-n,N}$, $B_2=0$, the proof follows from Theorem~\ref{thm:ErrorEstimate}.
\end{proof}
\begin{remark}
  Equation~\eqref{eq:EstRIP} gives an estimate on the solution of
  the $\ell_1$-mini\-miz\-ation problem. If we assume that $\sigma_k(x^N)_{1} \leq C_N
  \sigma_k(x^{\infty})_{1}$, which means that the best approximation of
 $x^N$ is as good as the best approximation of $x^{\infty}$, then
  Theorem~\ref{thm:ErrorEstimate} shows that the solution that we get
  from $\ell_1$-mini\-miz\-ation is as good as the best $k$-term approximation
  of $x^{\infty}$, where $x^{\infty}$ is the solution of original equation
  $A^{\infty} x^{\infty} = b^{\infty}$.
\end{remark}

\begin{remark}
  Equation~\eqref{eq:EstRIP} only gives a good bound, if we have
  \[
  (C^2_{2,k}-\frac{C^2_{1,k}}{k})\, \MutInc{B_1} \leq \frac{C_{1,k}^2}{k} (\mu_{\max}^2 - 1),
  \]
  where $\mu_{\max}$ is the largest singular value of $B_1$. Otherwise we
  may use the direct estimate $\norm{B_1 (x^N - \tilde{x})}{2}^2 \leq
  \norm{B_1}{2}^2 \norm{x^N - \tilde{x}}{2}^2$ and then apply
  Theorem~\ref{thm:CohenBestK1}.
\end{remark}

\section{Numerical Experiments}
\label{sec:numerics}

In this section we present some numerical examples.

\subsection{Example: 1D-Poisson Equation}
\label{sec:Example1D}

Let us first demonstrate that $\ell_1$-mini\-miz\-ation can successfully
obtain a sparse solution.  We consider the Poisson equation
\[
-u''=f \qquad\text{on}\qquad \Omega = (-1,1),
\]
with boundary conditions $u(-1) = 0 = u(1)$ and
\[
f(x) = 2 \cdot \alpha^3 \Big(
\frac{x+\tfrac{1}{2}}{1 + \alpha (x + \tfrac{1}{2})^2} +
\frac{x}{1+(\alpha \cdot x)^2} +
\frac{x - \tfrac{1}{4}}{1 + \alpha (x - \tfrac{1}{4})^2} +
\frac{x - \tfrac{1}{2}}{1 + \alpha^2 (x - \tfrac{1}{2})^2}
\Big),
\]
where $\alpha \define 100 \cdot \pi$. The exact solution of this problem is
\begin{align*}
  u(x) = & \arctan(\alpha (x + \tfrac{1}{2})) + \arctan(- \alpha (x + \tfrac{1}{2}))
  + \arctan(\alpha \cdot x) \\
  & + \arctan(- \alpha \cdot x) + \arctan(\alpha (x - \tfrac{1}{4})) +
  \arctan(- \alpha (x - \tfrac{1}{4})) \\
  & + \arctan(\alpha (x - \tfrac{1}{2})) + \arctan(- \alpha (x - \tfrac{1}{2})).
\end{align*}

\begin{figure}
  \includegraphics[width=0.45\textwidth]{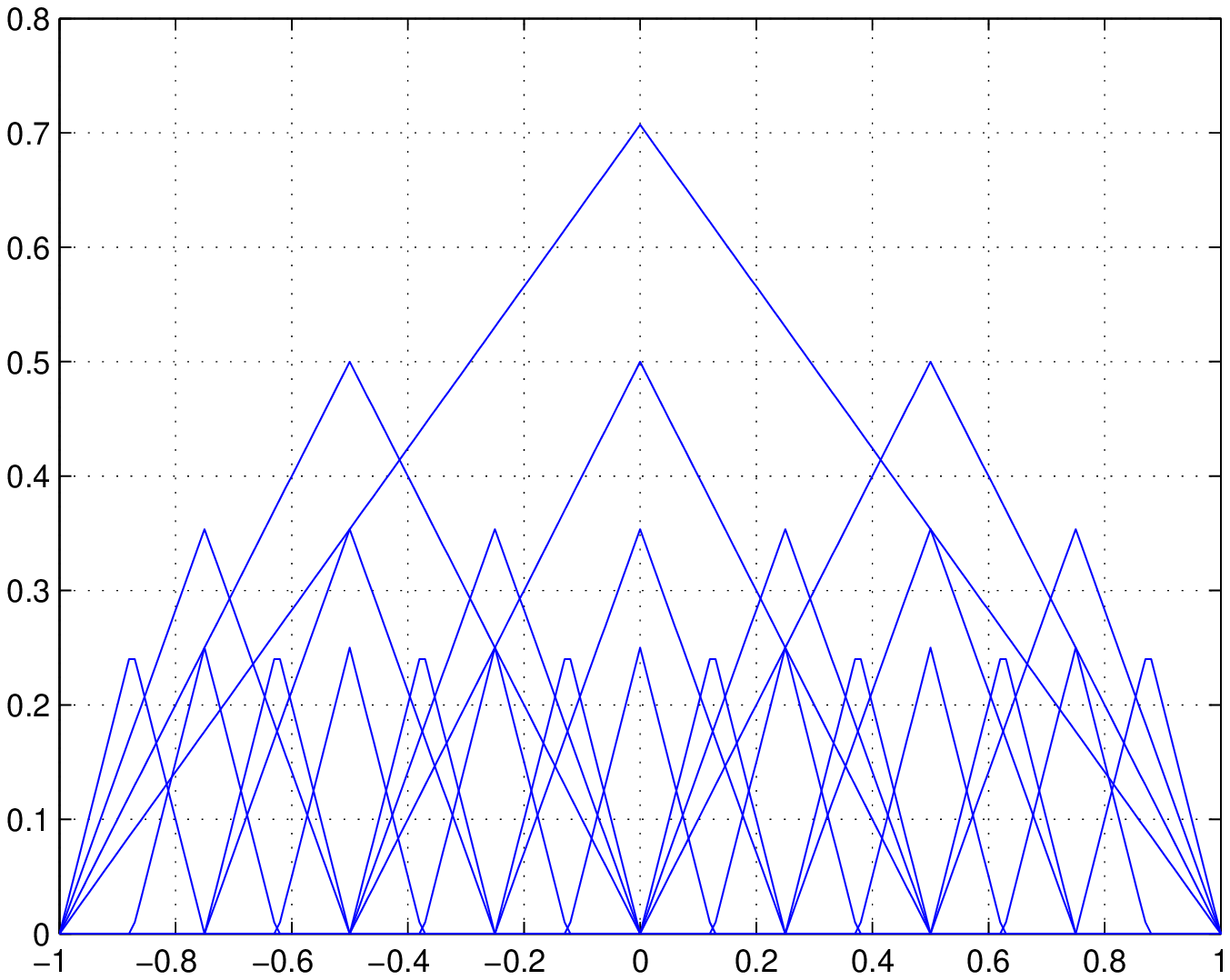}\hfill
  \includegraphics[width=0.45\textwidth]{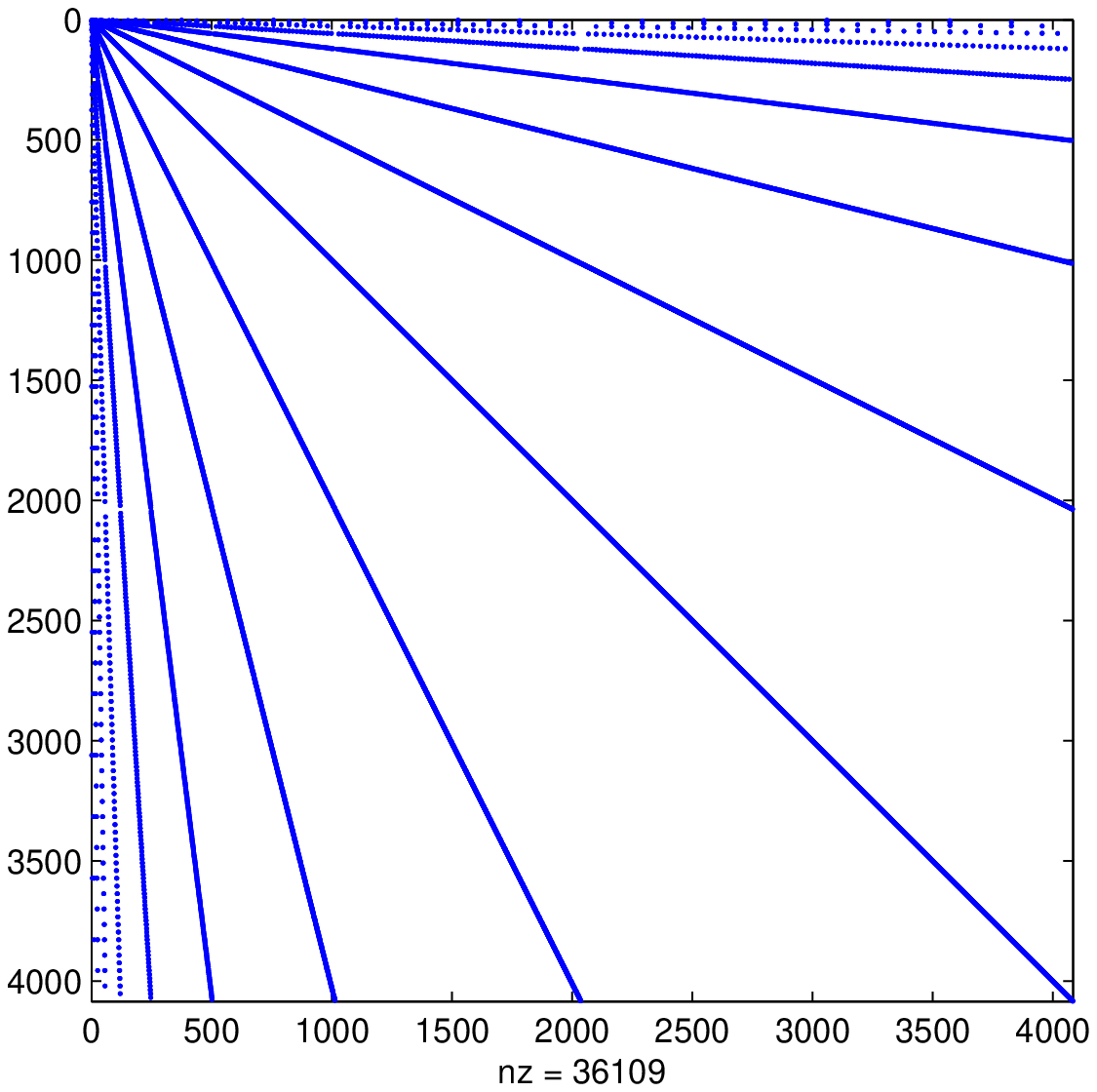}
  \caption{\emph{Left:} four levels of hat functions on the interval
    $[-1,1]$. \emph{Right:} nonzeros in the corresponding stiffness matrix
    with $11$ levels.}
  \label{fig:Multihat-11}
\end{figure}

We apply a finite element method~\cite{Bra07} and use classical \emph{shape
  functions}
\begin{equation}\label{MotherHat2}
  \phi(x) \define
  \begin{cases}
    1 - \abs{x} & \text{if }-1 \leq x \leq 1 \\
    0 & \text{otherwise.}
  \end{cases}
\end{equation}
The different refinement levels are given by
\[
\phi_{k,\ell}(x) \define 2^{-k/2}\, \phi\big(2^{k-1}\,(x+1)-\ell\big),
\ k \in \N,\; \ell = 1, \dots,2^{k}-1.
\]
Here, the scaling factor $2^{-k/2}$ is used to make the diagonal elements
of the stiffness matrix equal to $1$. We then have
\[
\phi_{k,\ell}(x) = 2^{-k/2}
\begin{cases}
  1 - \abs{2^{k-1}\, (x+1) - \ell} & \text{if }
  -1 + \frac{\ell-1}{2^{k-1}} \leq x \leq -1 + \frac{\ell+1}{2^{k-1}} \\
  0 & \text{otherwise,}
\end{cases}
\]
see the left part of Figure~\ref{fig:Multihat-11} for an illustration.
On the $N$th level, we then have
\[
\sum_{k=1}^N (2^k-1) = 2^{N+1}-(N+2)
\]
generating functions and we obtain a stiffness matrix that has a sparsity
as depicted on the right part of Figure~\ref{fig:Multihat-11}.

For this problem we have  numerically computed the mutual incoherence and the
restricted isometry constant. The numerical results indicate
that for $A^{N-n,N}$ as in (\ref{eq:AN-n}) we have
\[
\mathcal{M}(A^{N-n,N}) = \sqrt{\tfrac{2}{3}},
\]
independently of the size of the matrices, and that $\delta_1 = 0$,
$\delta_2 = 0.8165$, and for all $k > 1$, we have $\delta_k > 1$.

\begin{figure}
  \includegraphics[width=0.45\textwidth]{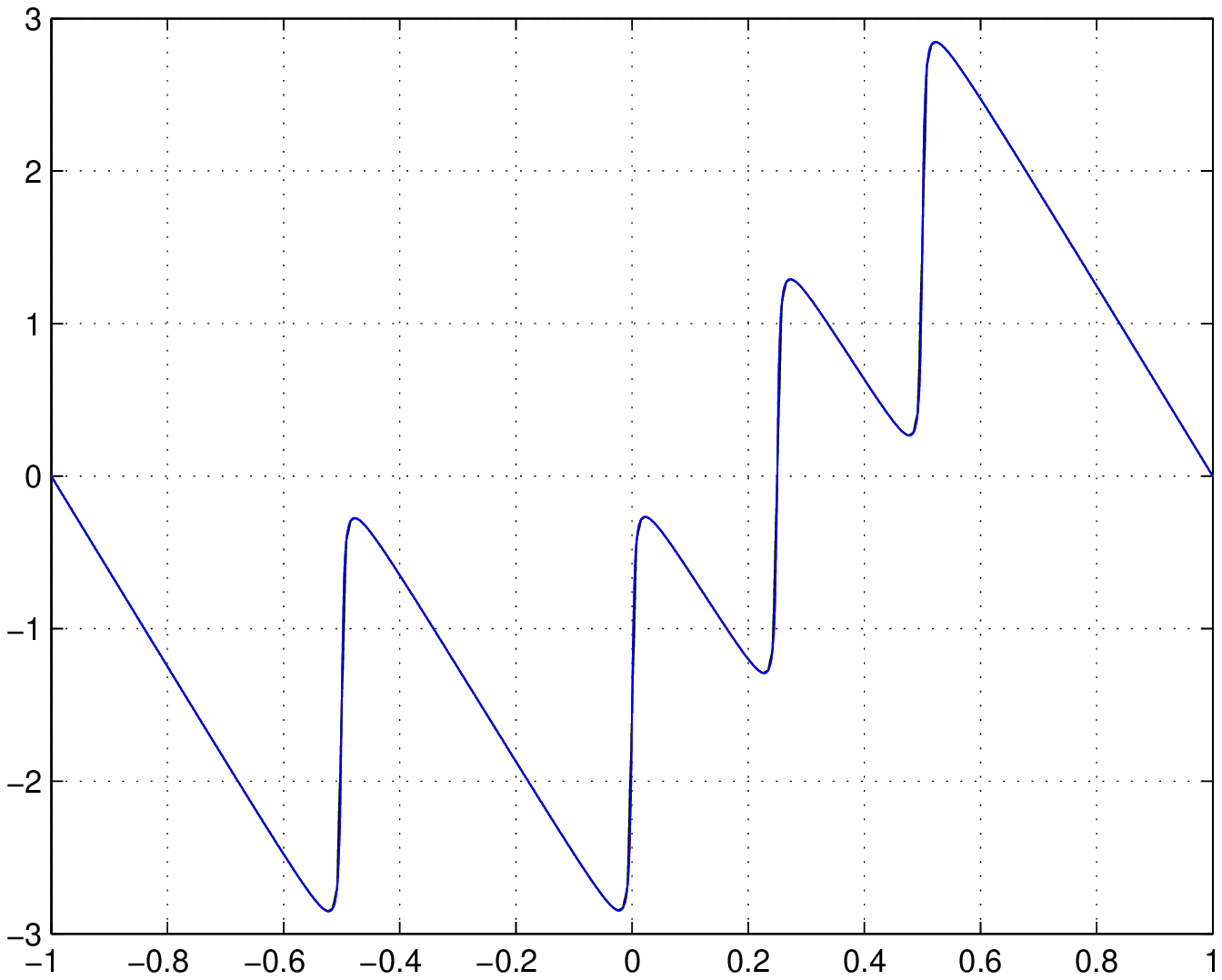}\hfill
  \includegraphics[width=0.45\textwidth]{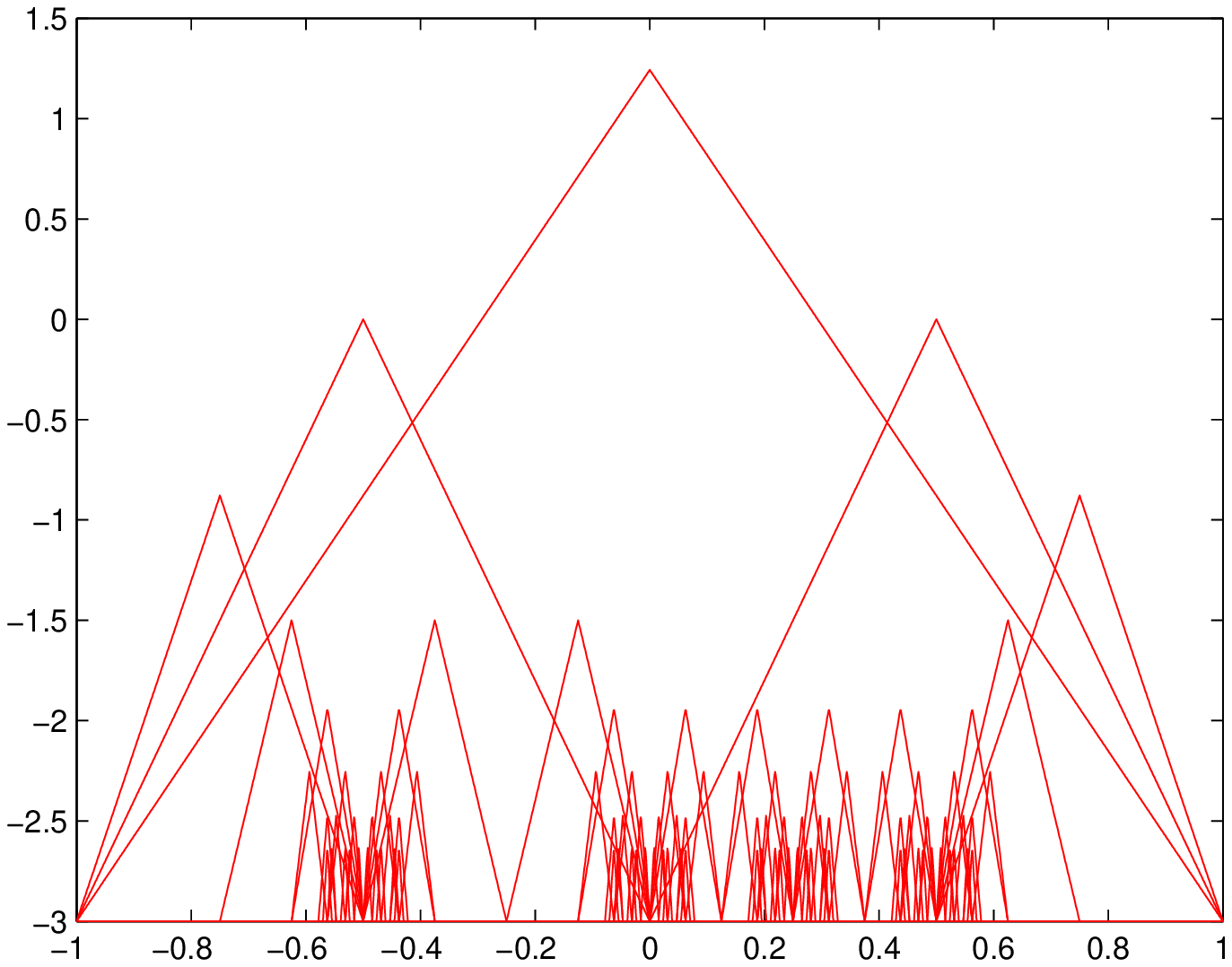}
  \caption{\emph{Left:} exact solution and approximate solution of the
    example in Section~\ref{sec:Example1D} at level 8.  \emph{Right:} hat
    functions that are selected by the $\ell_1$-mini\-miz\-ation problem.}
  \label{fig:Multihat-11result}
\end{figure}

On level $8$ we used the matrix $A^{N-n,N}$ of size $255 \times 502$. The
least squares solution of $A^N x = b$ has \emph{495} nonzeros, while the
minimum $\ell_1$-solution only had \emph{57} nonzeros.  The left part of
Figure~\ref{fig:Multihat-11result} depicts the exact solution and the
approximate solution at level~$8$.  There is no obvious difference and the
relative error in $\ell_2$-norm is $0.0627$.  The right part of
Figure~\ref{fig:Multihat-11result} shows that our method refines properly
at points with large gradients.

\subsection{Application of Algorithm~\ref{alg:IRBP} to a 1D-Poisson Equation}
\label{Iter:Exam1}

To illustrate the behavior of Algorithm~\ref{alg:IRBP}, we consider the
Poisson equation
\begin{equation}\label{Poissoneq4}
  \begin{array}{l l}
    Lu=-u''=2\frac{(100\pi)^3 x}{1+(100\pi x)^2}+2\frac{(100\pi)^3
      (x-0.5)}{1+(100\pi (x-0.5))^2}, \quad x\in \Omega=(-1,1),\\
    (u(-1),u(1))=(0,0).
  \end{array}
\end{equation}
The exact solution of this problem is
\begin{align*}
  u = & \arctan(100\, \pi \cdot x) + \arctan(-100\, \pi)\cdot x \\
  & + \arctan(100\, \pi (x - \tfrac{1}{2}))
  + \arctan(- 100\, \pi) (x - \tfrac{1}{2}).
\end{align*}
We applied four refinement steps of Algorithm~\ref{alg:IRBP} starting from level~4. In turns out that starting from level 3, a
straightforward refinement process does not work, since some of the
singularities are lost; this is a point where our algorithm would need to
backtrack, see Section~\ref{sec:Algorithm}.

For the starting level $4$, the size of matrix $A_{11}$ in~\eqref{eq:Partition} is $2^{r}-(r+1)=11$ for $r=4$. The size of $A_{21}$ is $(2^r-1) \times (2^{r}-(r+1))=15 \times 11$. The size of $A_{22}$ is  $15\times 15$.

\begin{figure}
  \centering
  \includegraphics[angle=0, width=0.35\textwidth]{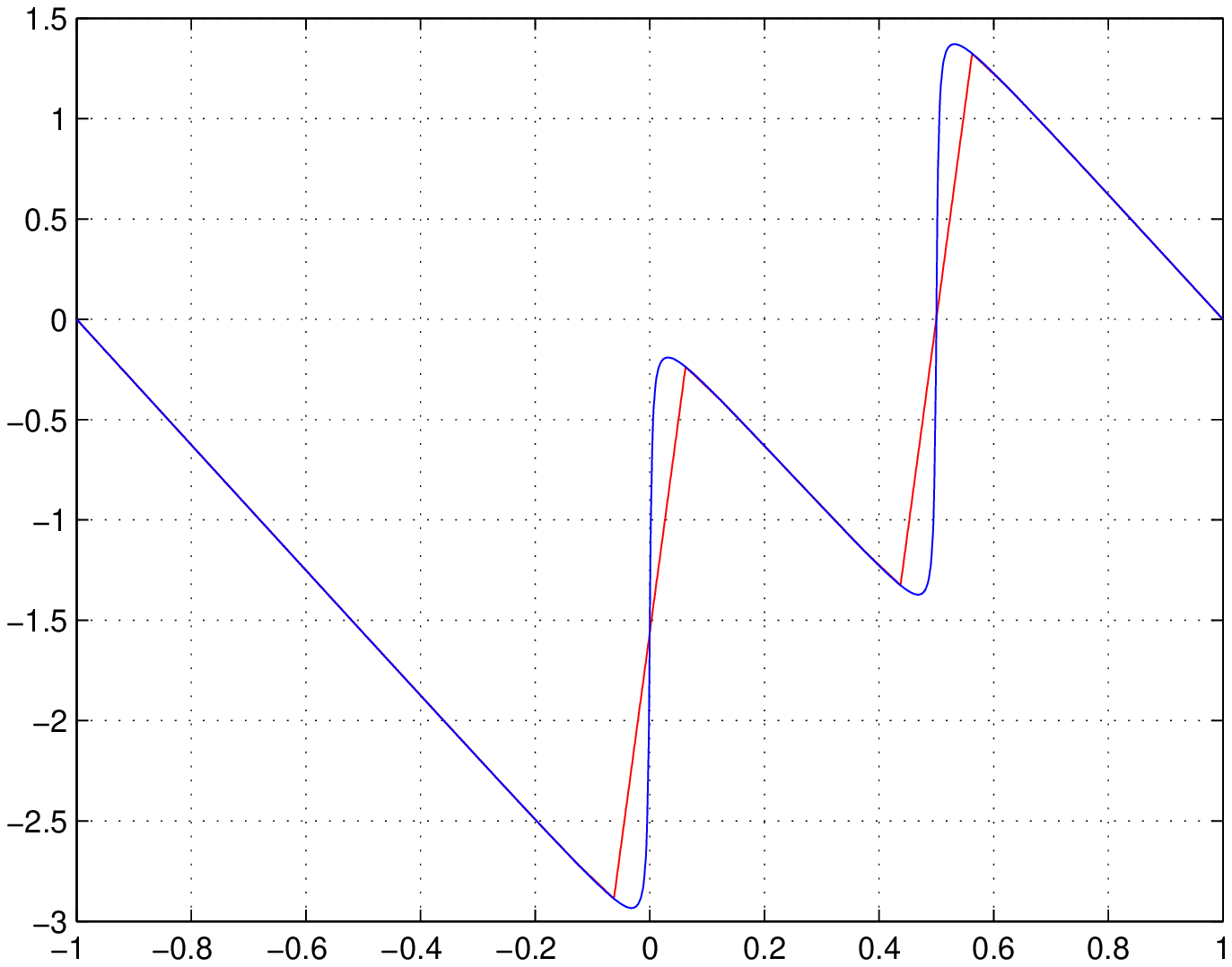}\qquad
  \includegraphics[angle=0, width=0.35\textwidth]{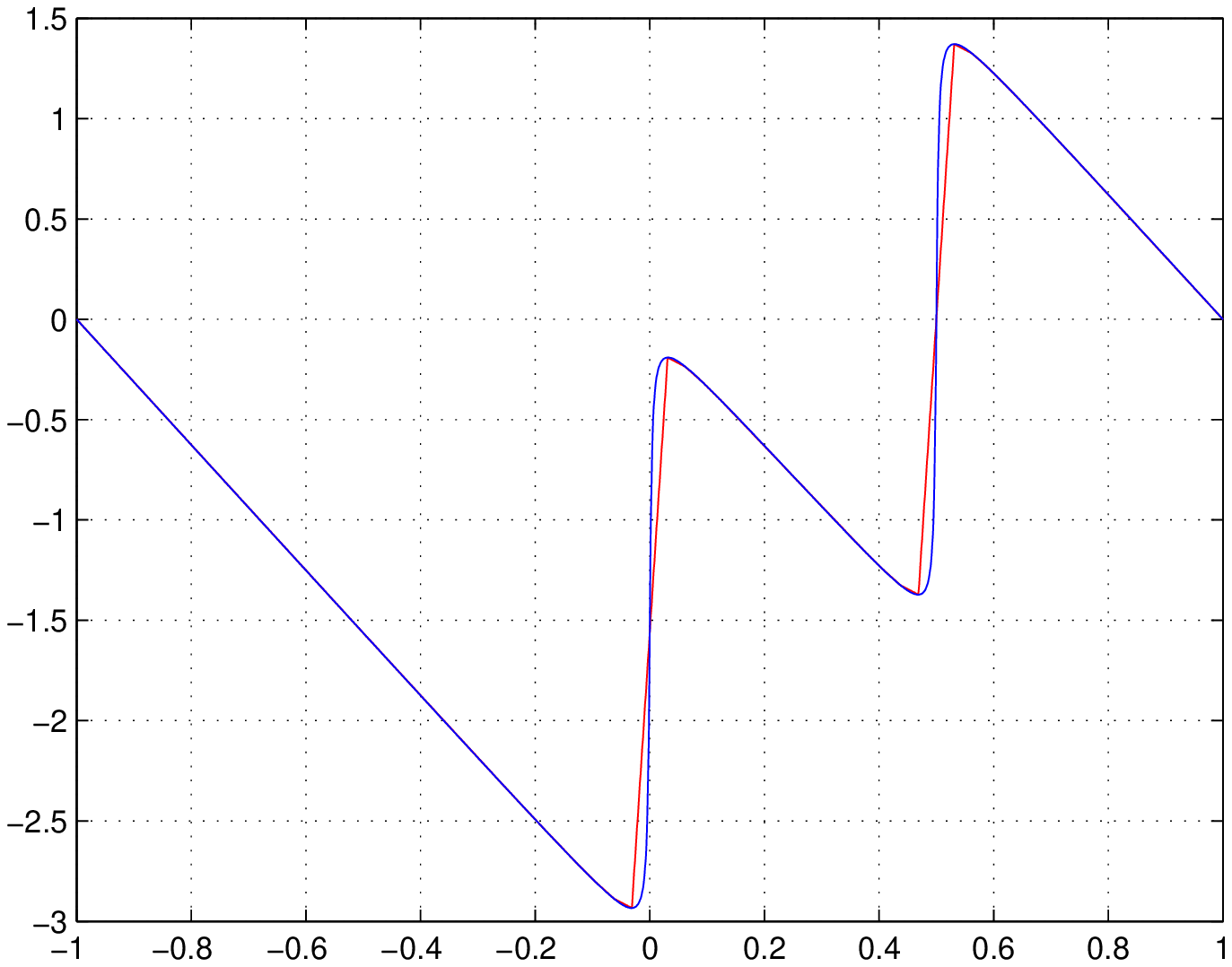}\\[1ex]
  \includegraphics[angle=0, width=0.35\textwidth]{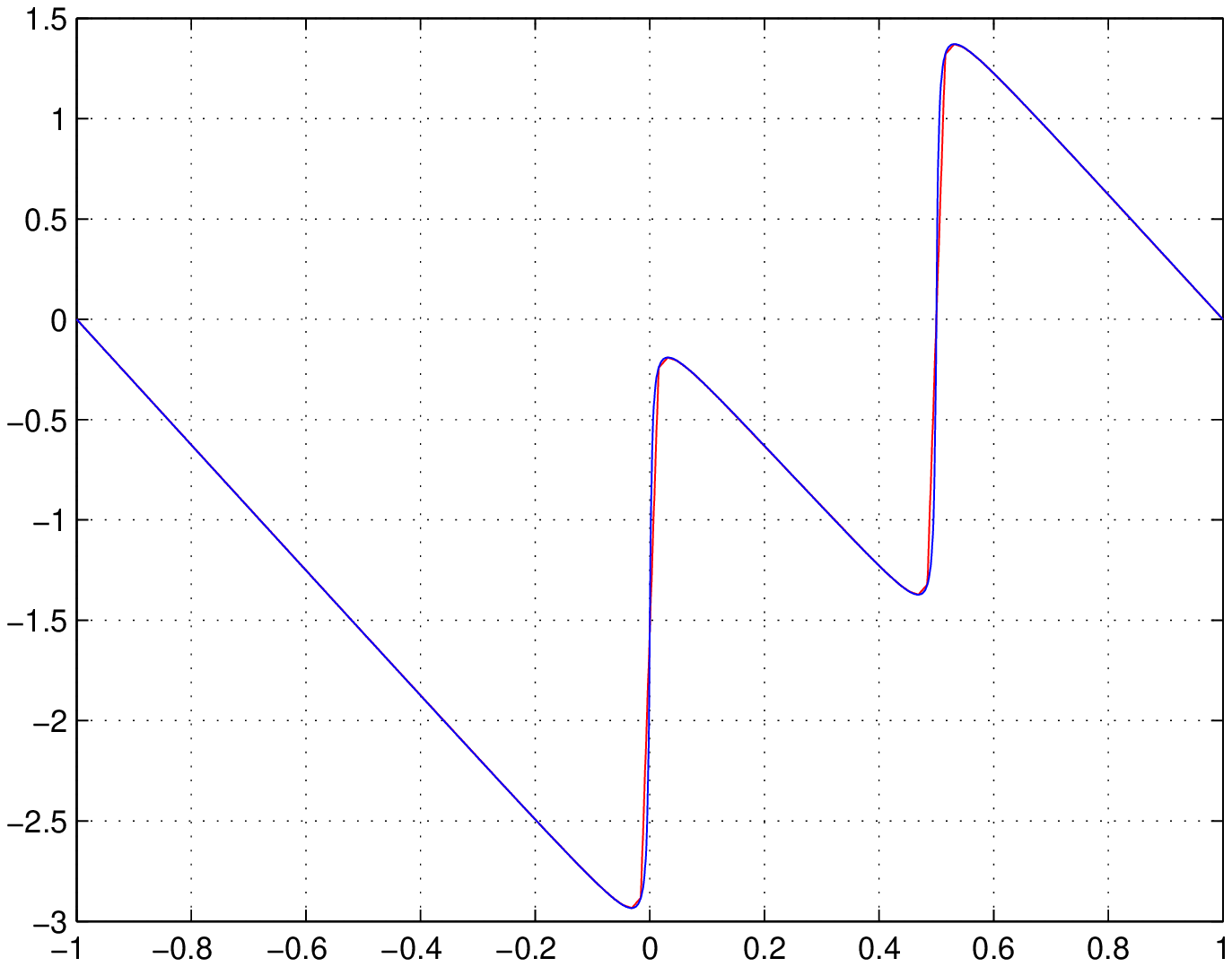}\qquad
  \includegraphics[angle=0, width=0.35\textwidth]{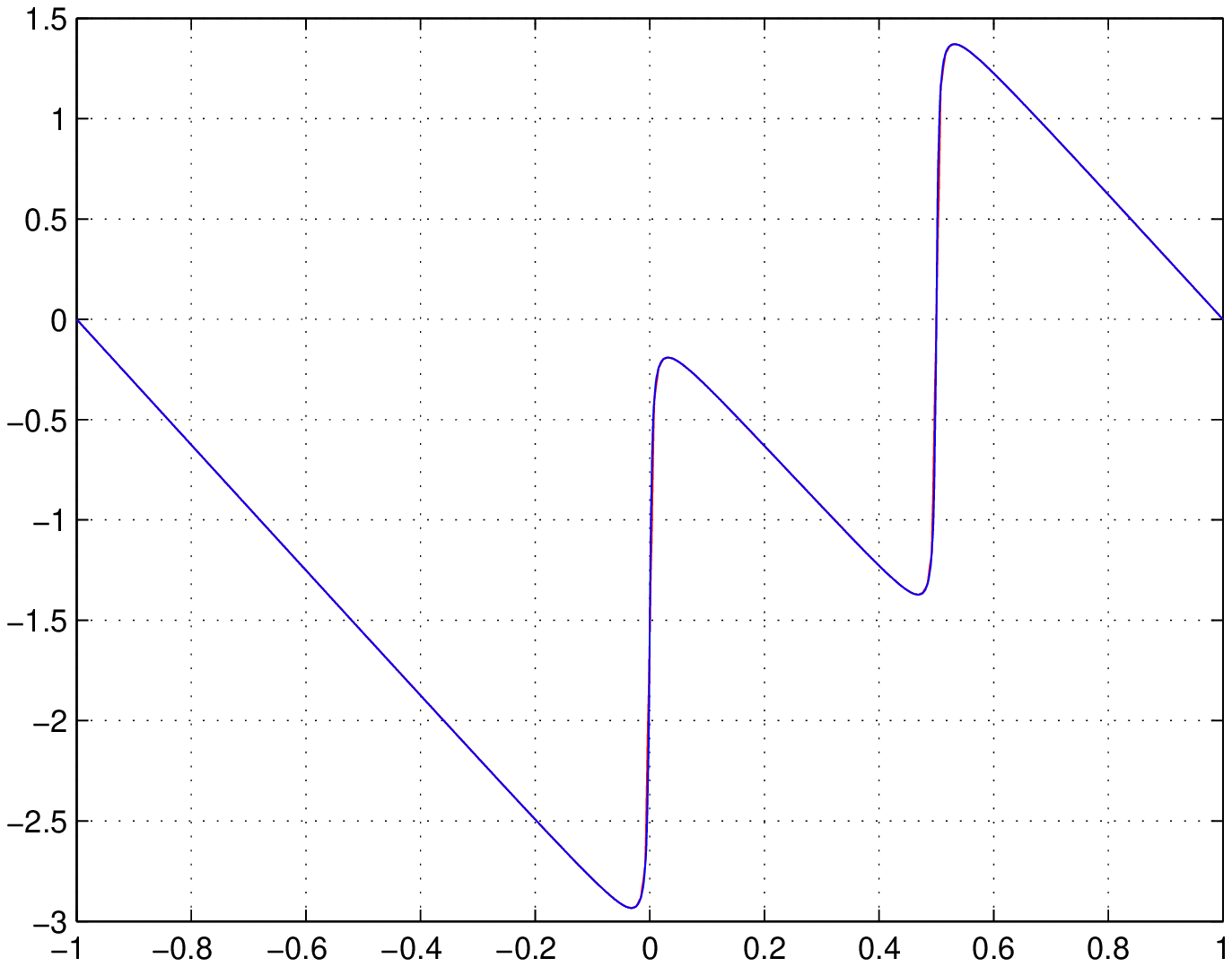}
  \caption{Exact and approximate solution for Example~\ref{Iter:Exam1} using Algorithm~\ref{alg:IRBP}.}
  \label{IterMulhatres1}
\end{figure}

\begin{table}[tb]
  \caption{Results of Algorithm~\ref{alg:IRBP}
    for Example \ref{Iter:Exam1}. The first column is the refinement step,
    column 2 gives the size of the matrix used for the $\ell_1$-mini\-miz\-ation problem,
    column 3 gives the $\ell_0$-norm of the $\ell_1$-minimal solution,
    column 4 gives the size of the matrix used for the classical FEM,
    column 5 gives the FEM solution of the problem at different levels,
    and column 6 gives the ratio of $\ell_0$-norms of the $\ell_1$-solution and the FEM solution.}

  \label{Iter:BCH0}
  \begin{tabular*}{\textwidth}{@{\extracolsep{\fill}}rrrrrr@{}}\toprule
    step  & size of $\ell_1$-matrix  &  $\norm{z}{0}$   &  size of FEM matrix & $\norm{x}{0}$
      & $\norm{z}{0}/\norm{x}{0}$ 
\\\midrule
    1 & 15 $\times$ 26   &  13 & 15 $\times$ 15 &  15 & 0.867 \\
    2 & 29 $\times$ 42   &  23 & 31 $\times$ 31 &  31 &0.742 \\
    3 & 59 $\times$ 72   &  41 & 63 $\times$ 63 & 63 & 0.651 \\
    4 & 113 $\times$ 126 &  67 & 127 $\times$ 127 & 127 & 0.528 \\
    5 & 191 $\times$ 204 & 103 & 255 $\times$ 255 & 255 & 0.404 \\\bottomrule
  \end{tabular*}
\end{table}

At each step we determine the new support using $\ell_1$-mini\-miz\-ation
and then refine these nodes and all necessary higher level nodes according
to Algorithm~\ref{alg:IRBP}, see Figure~\ref{IterMulhatres1}.  In
Table~\ref{Iter:BCH0} we present the results of four refinement steps of
Algorithm~\ref{alg:IRBP}.

\subsection{Application of Algorithm~\ref{alg:IRBP} to a 2D-Poisson
  Equation}
\label{Iter:Exam2}

As a second example for Algorithm~\ref{alg:IRBP}, we consider the Poisson
equation
\[
-\frac{\partial^2 u}{\partial x^2}-\frac{\partial^2 u}{\partial y^2}=f(x,y)
\  \text{on} \  [0,1]\times [0,1]
\]
with $u(x,y)=0$ on the boundary of $[0,1]\times [0,1]$ where $f(x,y)=-20x(x-1)-20y(y-1)$.
The original solution is $u(x,y)=10xy(x-1)(y-1)$.

\begin{figure}
  \centering
  \includegraphics[angle=0, width=0.3\textwidth]{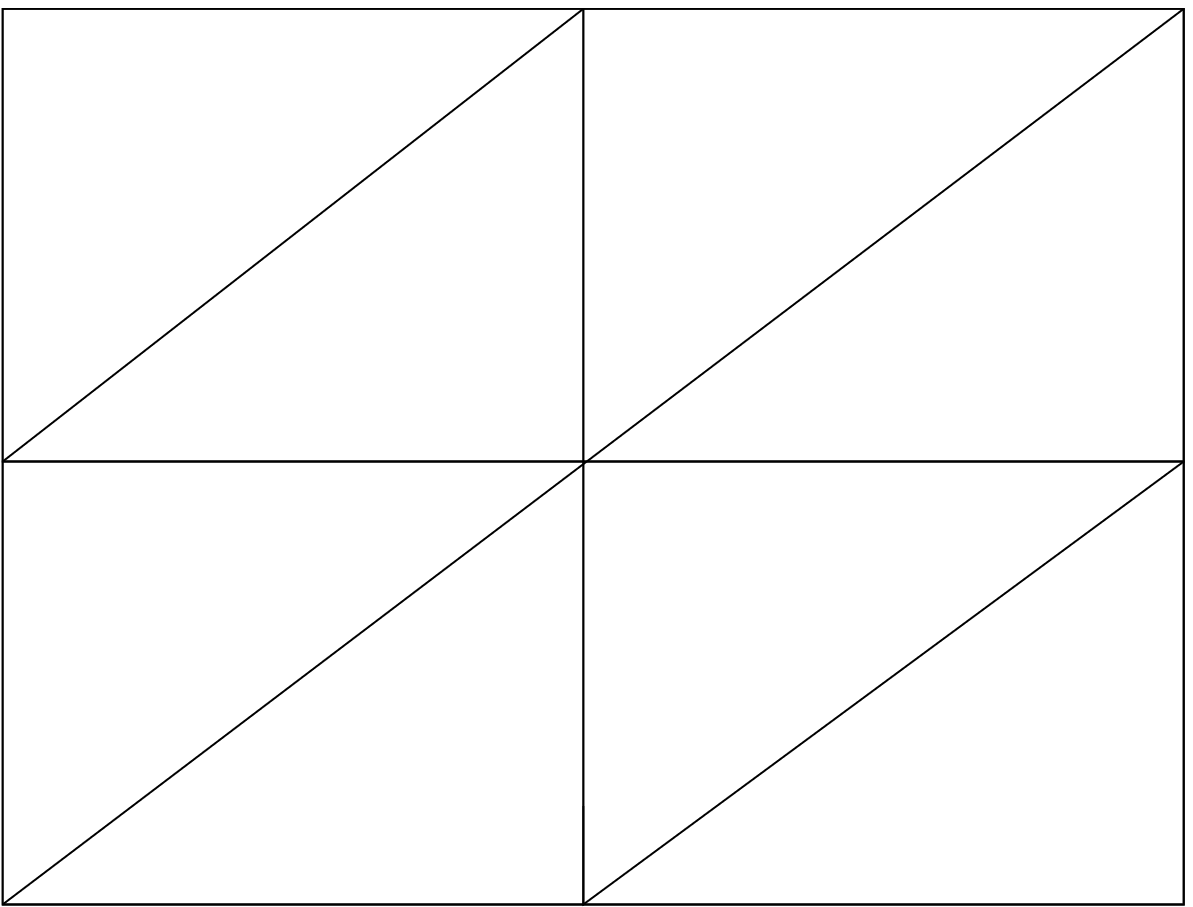}\qquad
  \includegraphics[angle=0, width=0.3\textwidth]{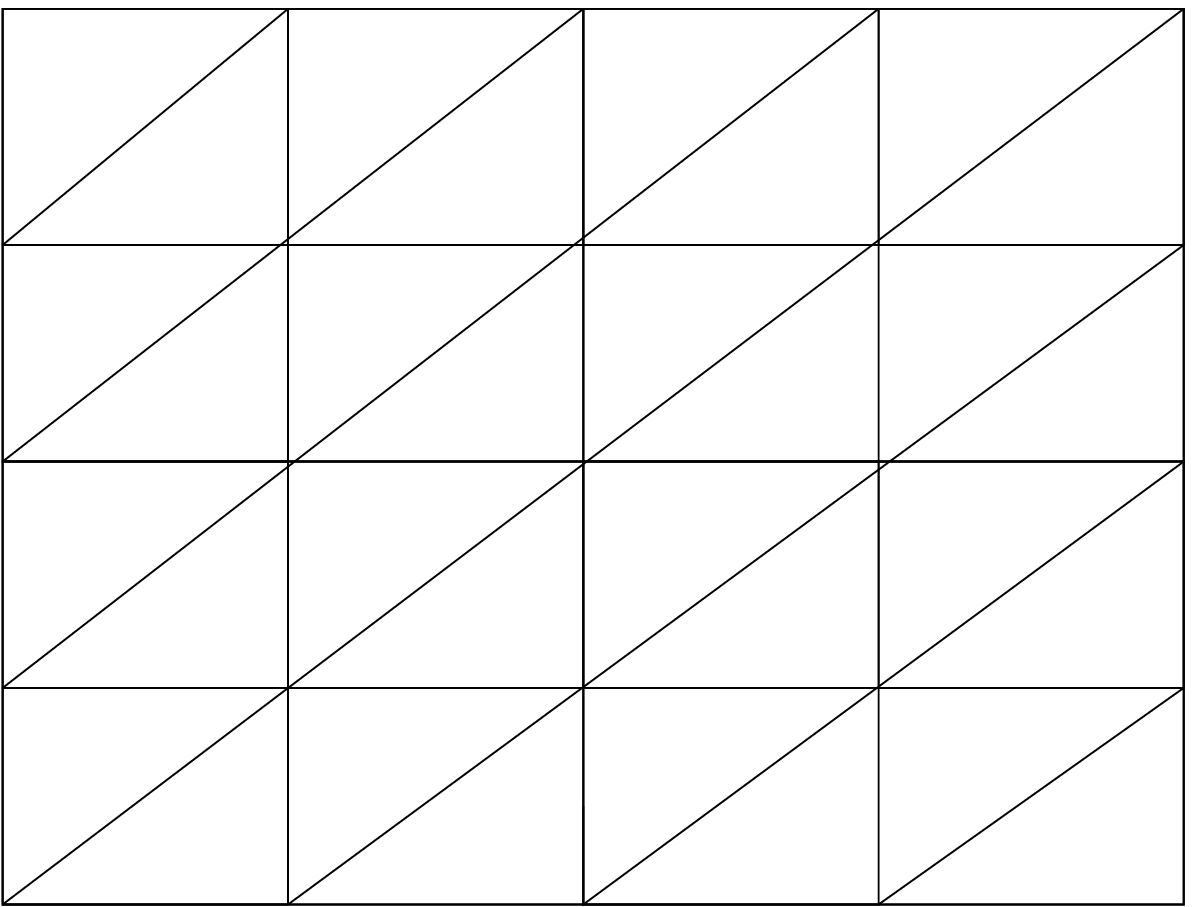}
  \caption{Triangulation on first and second level}
  \label{fig:Triangulation}
\end{figure}

\begin{figure}
  \centering
  \includegraphics[angle=0, width=0.6\textwidth]{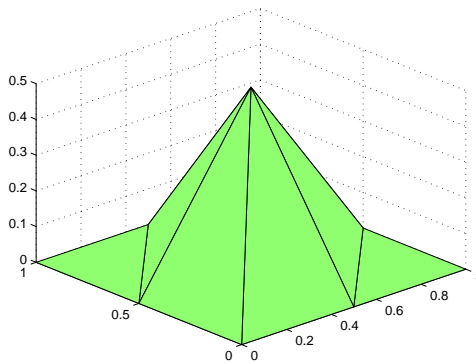}
  \caption{Basis function in the first level}
  \label{fig:Triangulation1}
\end{figure}

We use piecewise linear generating functions of the form
$u_i(x,y)=a_ix+b_iy+c_i$ on each triangle. Figure~\ref{fig:Triangulation}
depicts the refinement step from the first level (left) to the second level
(right). The basis function on level~1 is plotted in
Figure~\ref{fig:Triangulation1}.

Based on the triangulation on level $j$ and the triangles
$\{\Delta_{j,k_i}\}_{i=1}^6$, we obtain the generating functions:
\[
  \phi_{j,(a,b)}(x,y) = \frac{1}{2}
  \begin{cases}
      2^j(x-a)+1 & (x,y)\in \Delta_{j,k_1} \\
      2^j(x-a)-2^j(y-b)+1 & (x,y)\in\Delta_{j,k_2} \\
      -2^j(y-b)+1 & (x,y)\in\Delta_{j,k_3} \\
      -2^j(x-a)+1 & (x,y)\in\Delta_{j,k_4} \\
      -2^j(x-a)+2^j(y-b)+1 & (x,y)\in\Delta_{j,k_5} \\
      2^j(y-b)+1 & (x,y)\in\Delta_{j,k_6} \\
    \end{cases}
\]
where $(a,b)$ is the center of the basis function. In general on level $n$ we have $(2^n-1)^2$ basis functions.

\begin{figure}
  \centering
  \includegraphics[angle=0, width=0.45\textwidth]{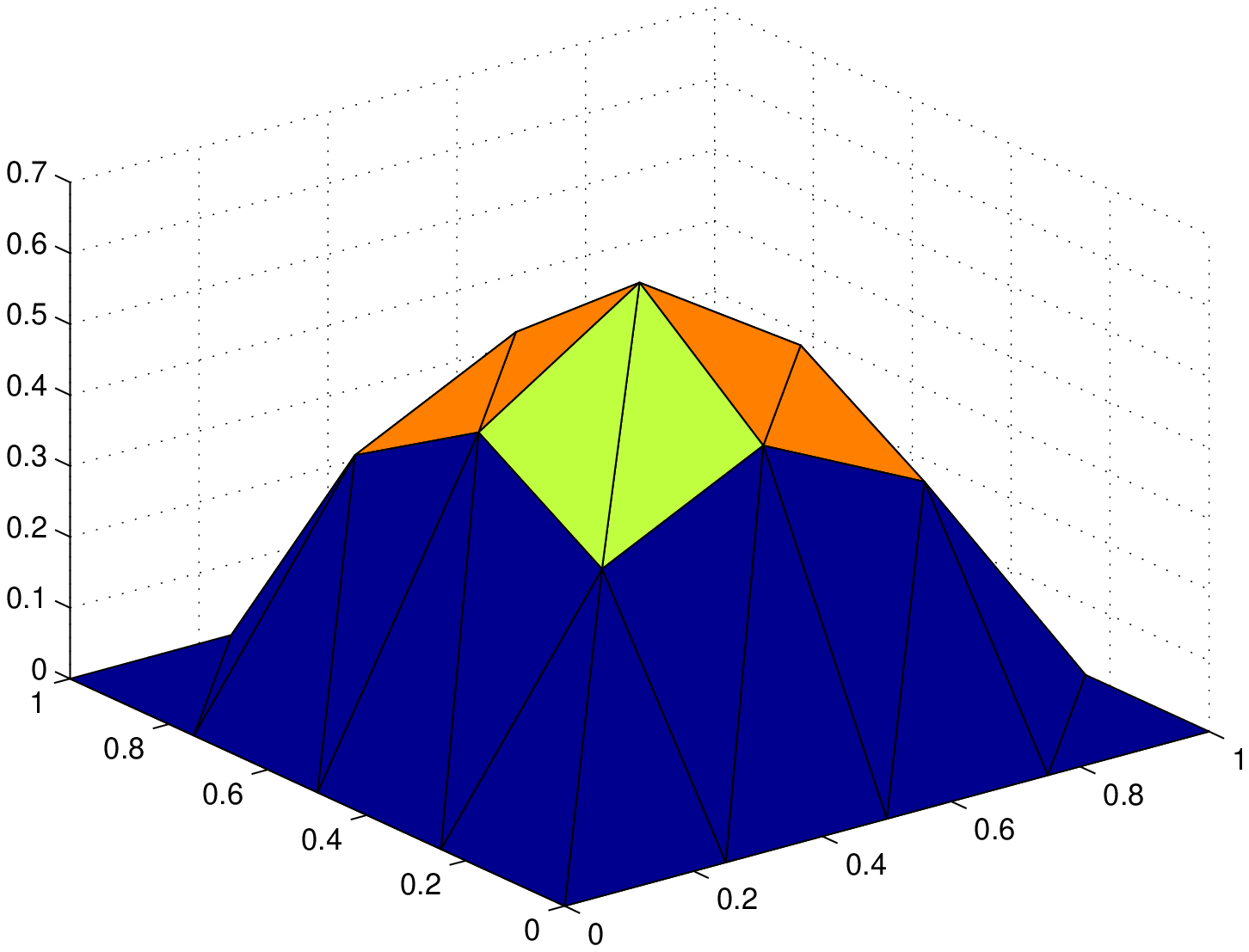}
  \includegraphics[angle=0, width=0.45\textwidth]{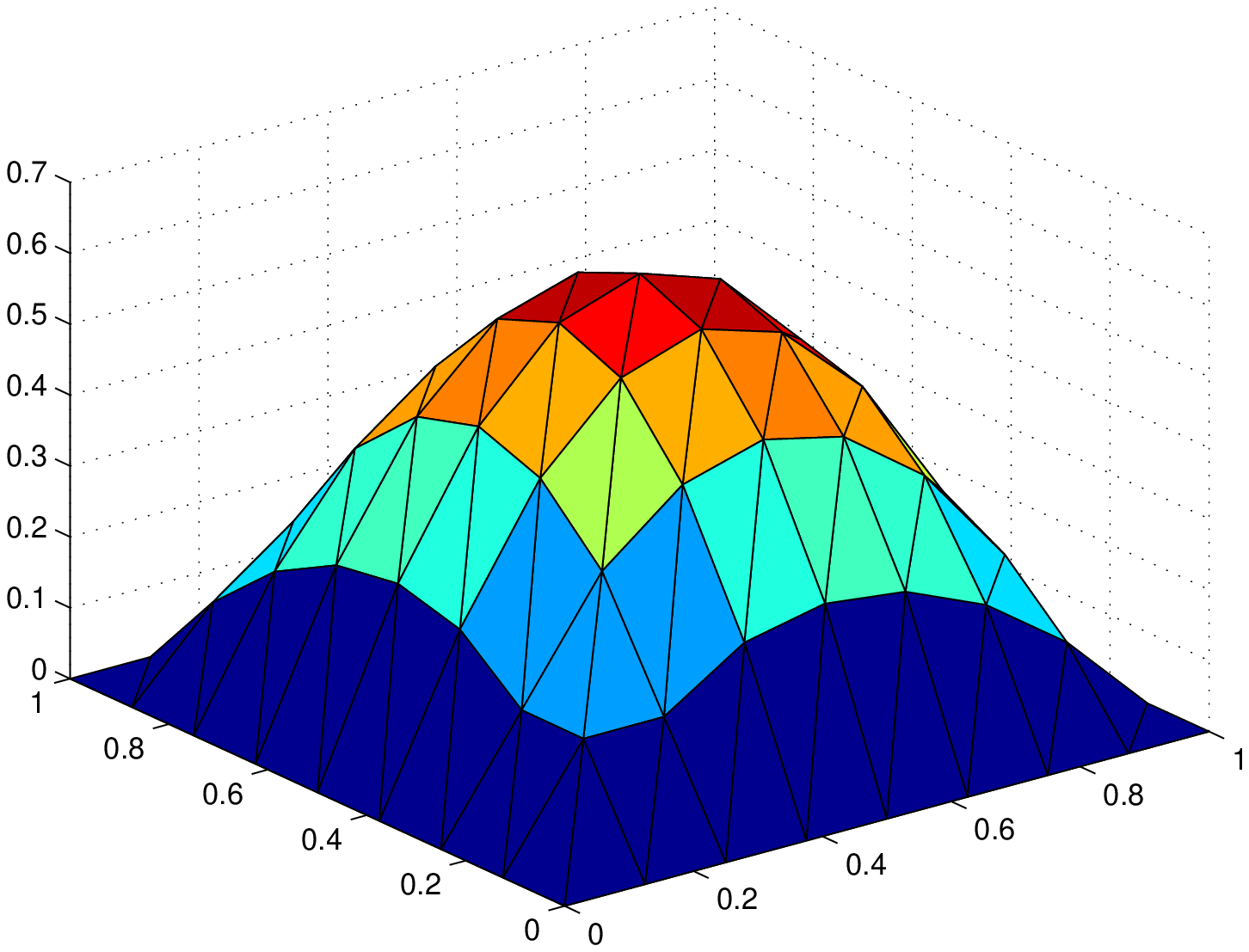}
  \includegraphics[angle=0, width=0.45\textwidth]{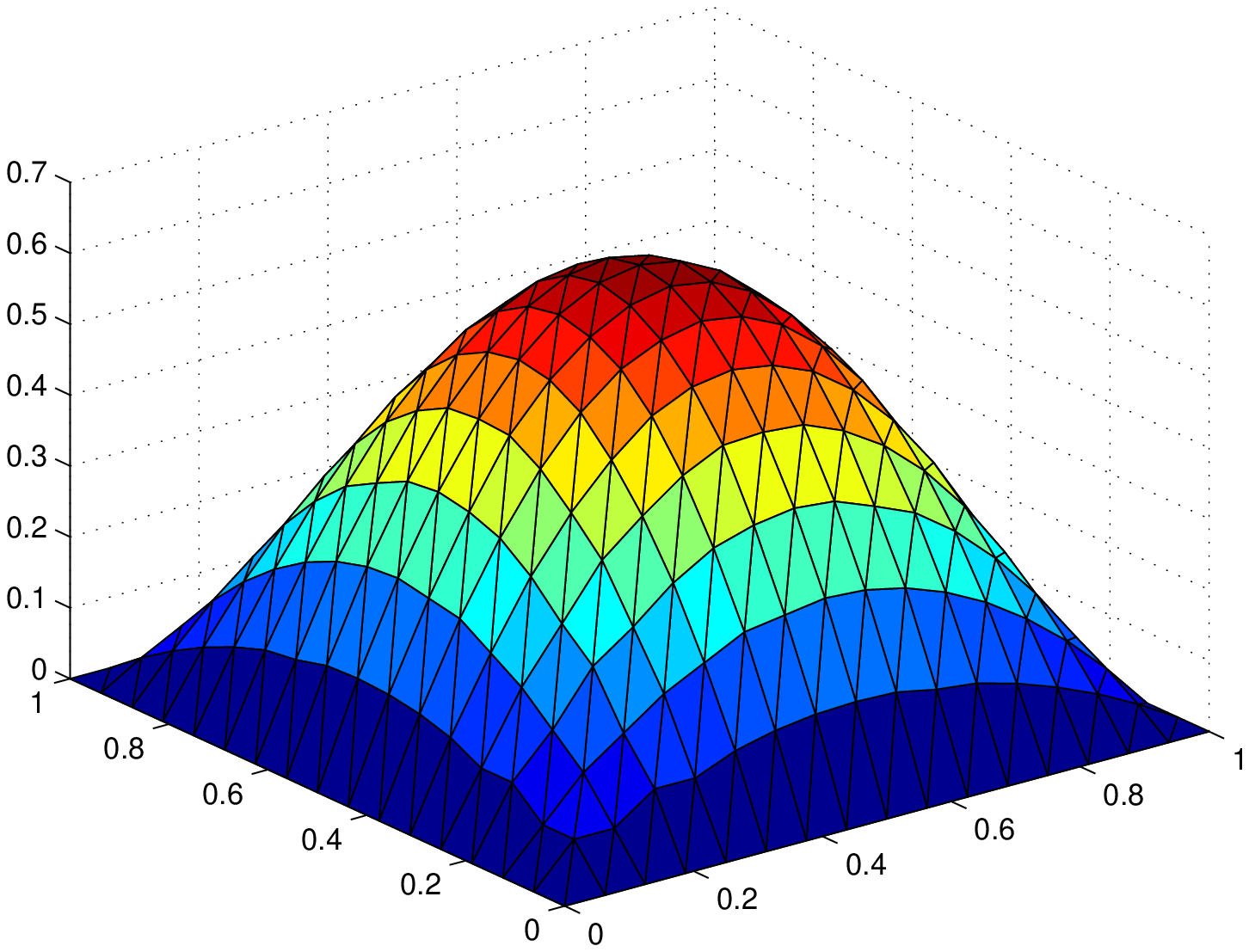}
  \includegraphics[angle=0, width=0.45\textwidth]{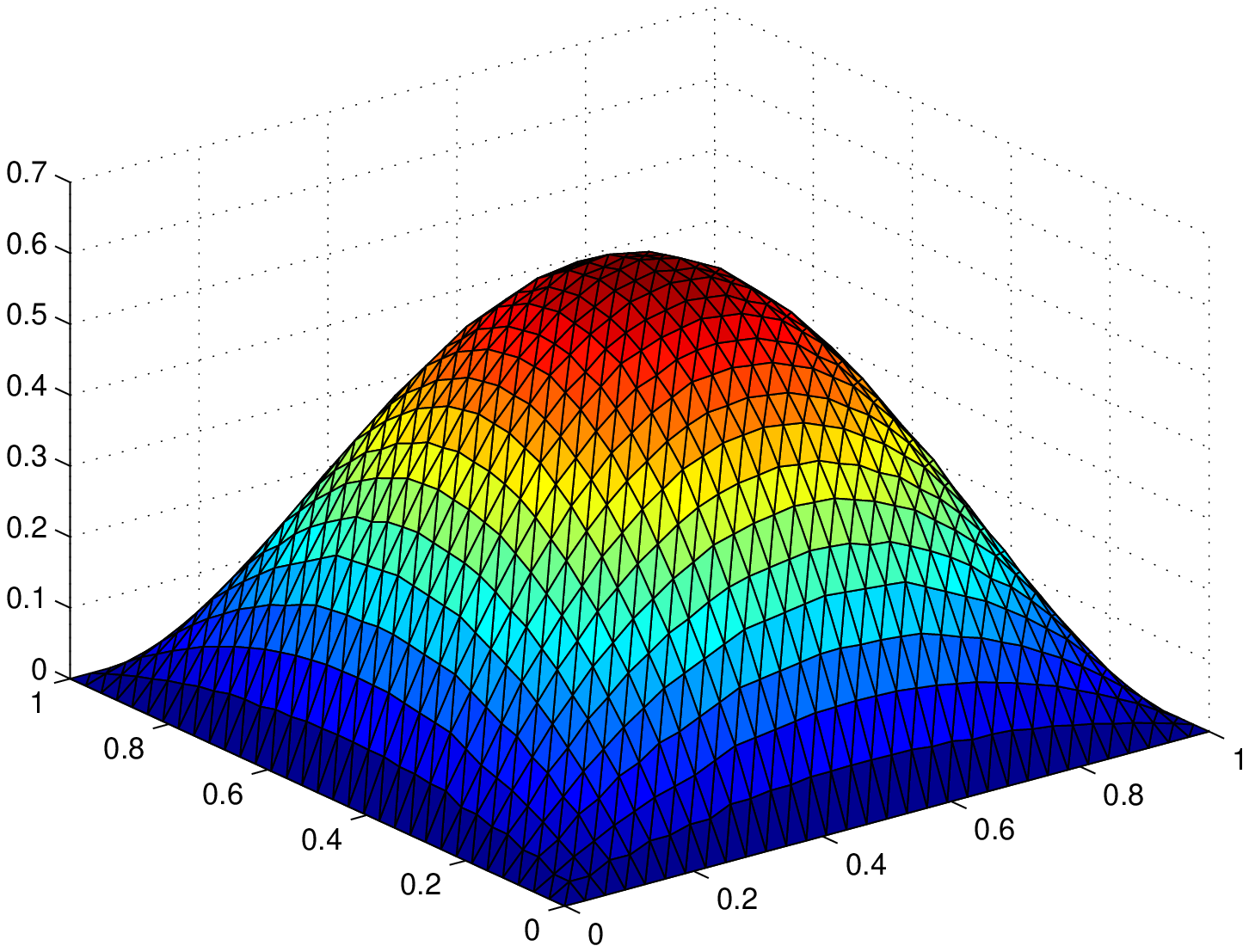}
  \includegraphics[angle=0, width=0.45\textwidth]{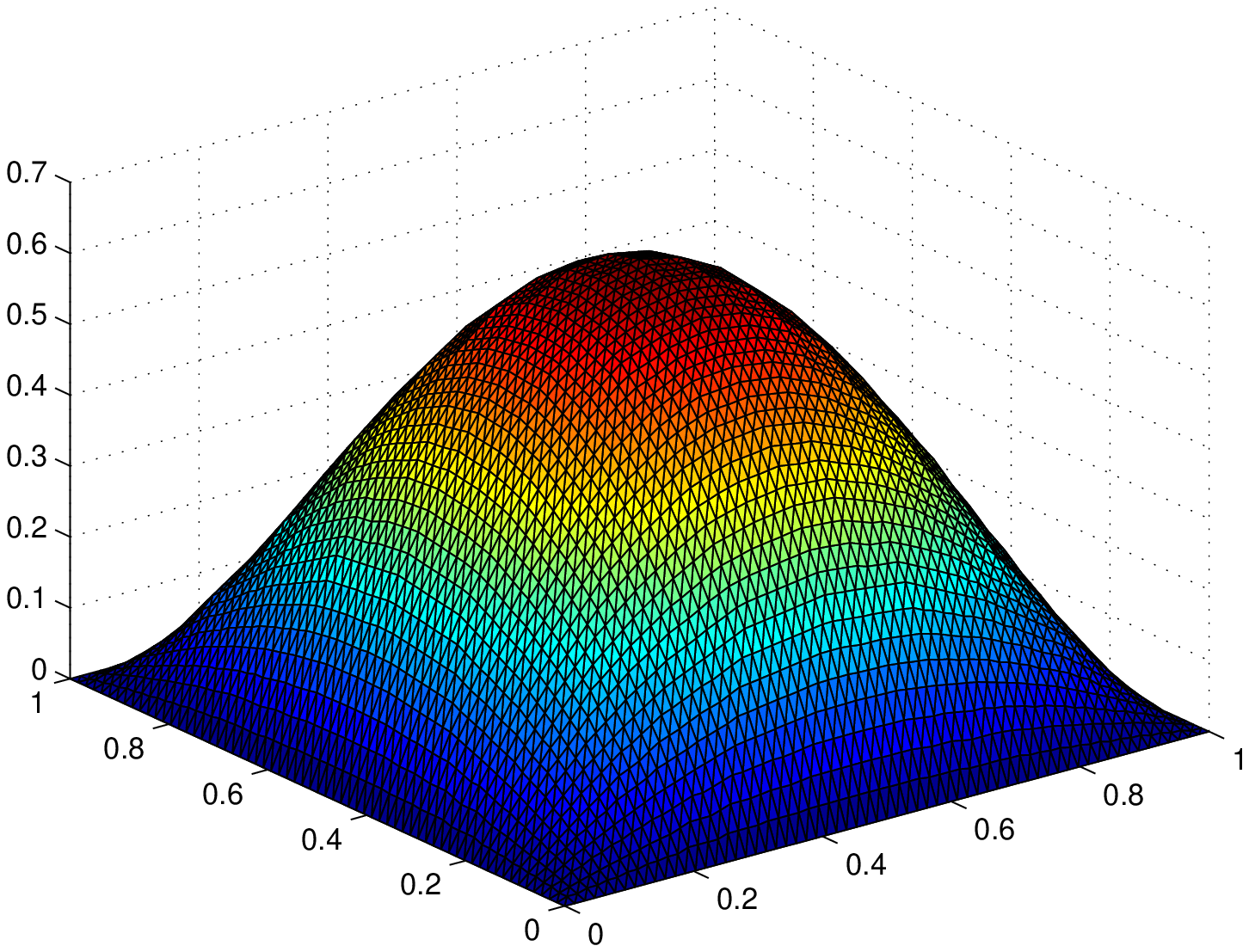}
  \includegraphics[angle=0, width=0.45\textwidth]{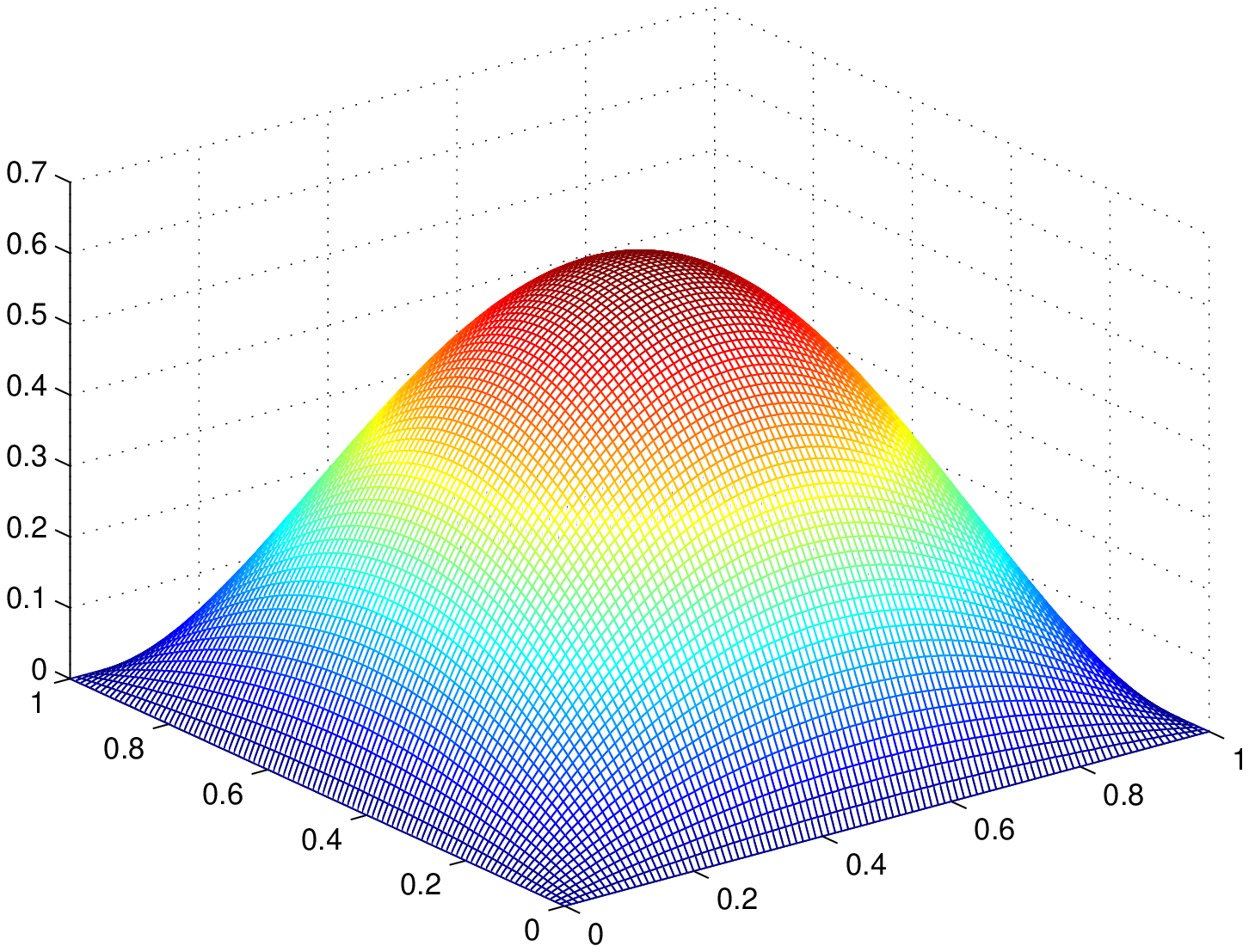}
  \caption{Approximate solutions obtained by Algorithm~\ref{alg:IRBP} on
    levels $2,3,4,5,6$ and exact solution (last) for Example~\ref{Iter:Exam2}.}
  \label{fig:IterMulhatres2}
\end{figure}

\begin{table}
  \caption{Results of Algorithm~\ref{alg:IRBP} for Example~\ref{Iter:Exam2}.
    The first column is the refinement step,
    column 2 gives the size of the matrix used for the $\ell_1$-mini\-miz\-ation problem,
    column 3 gives the $\ell_0$-norm of the $\ell_1$-minimal solution,
    column 4 gives the size of the matrix used for the classical FEM,
    column 5 gives the FEM solution of the problem at different level,
   and column 6 gives the ratio of $\ell_0$-norms of the $\ell_1$-solution and the FEM solution.}
\begin{tabular*}{\textwidth}{@{\extracolsep{\fill}}rrrrrr@{}}\toprule
    step  & size of $\ell_1$-matrix  &  $\norm{z}{0}$   &  size of FEM matrix & $\norm{x}{0}$
      & $\norm{z}{0}/\norm{x}{0}$ 
\\\midrule
    1 &      9 $\times$   10   &     8 &    9 $\times$    9 &    9  &  0.889  \\
    2 &     43 $\times$   51   &    42 &   49 $\times$   49 &   49  &  0.857  \\
    3 &    237 $\times$  244   &   114 &  225 $\times$  225 &  225  &  0.507  \\
    4 &    816 $\times$  823   &   172 &  961 $\times$  961 &  961  &  0.179  \\
    5 &   1352 $\times$ 1359   &   190 & 3969 $\times$ 3969 & 3969  &  0.048  \\\bottomrule
  \end{tabular*}
\label{tab:2D}
\end{table}

We applied four refinement steps of Algorithm~\ref{alg:IRBP} starting as first step from level 2.  At each step we determine the new support using $\ell_1$-mini\-miz\-ation and then refine these nodes and all necessary higher level nodes according to
Algorithm~\ref{alg:IRBP} (see Figure~\ref{fig:IterMulhatres2}). In Table~\ref{tab:2D} we present the results of four refinement steps of
Algorithm~\ref{alg:IRBP}. For the starting level $2$, the size of matrix $A_{11}$ in~\eqref{eq:Partition} is $\sum_{i=1}^{r-1} (2^i-1)^2=1$ for $r=2$. The size of $A_{21}$ is $(2^{r}-1)^2 \times \sum_{i=1}^{r-1} (2^i-1)^2)=9 \times 1$. The size of $A_{22}$ is  $9\times 9$.

\section{Efficiency Estimation}

The algorithm as described in Section~\ref{sec:Algorithm} relies on the
solution of a linear program (LP) in each iteration. Thus, for it to be
successful in practice the savings in the size of the matrices have to be
large enough to compensate for the higher solution times of LPs compared to
the standard finite element methods. Clearly, this depends on a number of
factors that are hard to predict: the PDE, the data, the required accuracy,
the concrete implementation, the usage of geometric refinement processes
etc.  This results in different sizes of the refinement tree and different
matrices with varying degrees of sparsity. We will nevertheless try to get
a rough idea of the efficiency. Since our method can be seen as a generic
way of controlling the refinement process, we compare it against the
classical finite element method without refinement. We use the examples of
Section~\ref{sec:numerics} as a guideline.

Note that we restrict our attention to the case of solving LPs, although
there are different methods for compressed sensing that yield similar
results as the LP-based methods; for instance, orthogonal matching pursuit
might be applied, see Tropp~\cite{Tro04}, Donoho, Elad, and
Temlyakov~\cite{DonET06}, Tropp and Gilbert~\cite{TroG07} and Needell and Tropp~\cite{NeeT08}. It remains
to be seen whether these methods are competitive for the application
discussed in this paper, especially with respect to sparse matrices.

Let us first consider worst case computing times. In the example in
Section~\ref{sec:Example1D}, the number of basis functions at level~$k$ is
$2^k-1$, which is also the size of the matrix $A^k$ in the equation
system~\eqref{eq:IterISystem}, which has to be solved by a ``classical''
finite element method. If we use a dense solver this takes $O(2^{3k})$ time
for each solution. In comparison, dense interior point algorithms for
linear programming require about $O(n^{3.5} L)$ time for an LP of
dimension~$n$, where $L$ is the encoding size of the LP. If the LP is
dense, the encoding length includes at least one bit for each entry of the
matrix and thus is of size at least $nm$, where~$m$ is the number of
constraints. In our case, we have $n = 2^{k+1}-(k+2)$ and $m = 2^k-1$, if
we would start our method at level~$k$. Hence, a very optimistic estimate
of the running time in the dense case would be $O(n^{3.5} mn) = O(2^{5.5 k})$.

\begin{table}[b]
  \caption{Running times for solving the LPs of the example in
    Section~\ref{sec:Example1D}; $m$ and $n$ denote the number of rows and
    columns of the constraint matrix, ``time'' refers to the running time
    in seconds, and $\norm{z}{0}$ gives the number of nonzeros in the
    solution.}
  \label{tab:RunningTimes}
  \begin{tabular*}{\textwidth}{@{\extracolsep{\fill}}rrrrr@{}}\toprule
    level & $m$  &   $n$ & time & $\norm{z}{0}$ \\\midrule
    7  &   127 &   247 &   0.01 & 7  \\
    8  &   255 &   502 &   0.03 & 8  \\
    9  &   511 &  1013 &   0.10 & 9  \\
    10 &  1023 &  2036 &   0.46 & 10 \\
    11 &  2047 &  4083 &   1.68 & 11 \\
    12 &  4095 &  8178 &   5.74 & 12 \\
    13 &  8191 & 16369 &  22.24 & 13 \\
    14 & 16383 & 32752 &  88.32 & 14 \\\bottomrule
  \end{tabular*}
\end{table}

Let us investigate the selection process of the compressed sensing approach
for this example, see Section~\ref{sec:Algorithm}. Assume that we are at
iteration~$k$ and we have the set $C^{k-1}$ of basis functions with $\ell_{k-1} =
\card{C^{k-1}}$. The refined set of basis functions $\hat{C}^k$ has size at
most $3 \, \ell_{k-1}$, because each basis function is subdivided into three new basis functions (some of the new basis functions might coincide). Assuming that we select at least a fraction of $\alpha \in
(0,1]$ among the basis functions of $\hat{C}^k$, the new iteration has at most
$\ell_{k} = \ell_{k-1} + \alpha\, 3\, \ell_{k-1} = (1 + 3\, \alpha) \ell_{k-1}$ basis
functions; compare Algorithm~\ref{alg:IRBP}. Therefore, at iteration~$k$,
we have at most $(1 + 3 \alpha)^k = 2^{k \log_2 (1 + 3 \alpha)}$ basis functions.
For the compressed sensing approach to be successful with respect to the finite element
(dense) worst case times above, we would need
\begin{align*}
& \log_2 (1 + 3 \alpha) < 3.5/5.5 \approx 0.64\\
\Leftrightarrow \quad & \alpha < \tfrac{1}{3}\, (2^{0.64} - 1) \approx 0.19.
\end{align*}
Hence, if we assume that the compressed sensing method yields a reduction
rate~$\alpha$ below approximately $0.19$, it should be effective, if we
assume worst case running times. In the results of Table~\ref{Iter:BCH0}, we have $\alpha$ around $0.5$. 


\begin{figure}
  \centering
  \includegraphics[angle=0, width=0.5\textwidth]{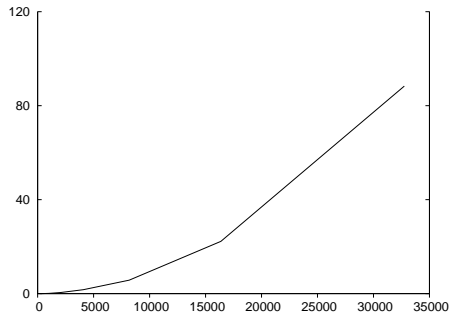}
  \caption{Solution times for solving the linear programs of
    Table~\ref{tab:RunningTimes} versus the number of columns of the
    matrix.}
  \label{fig:RunningTimes}
\end{figure}

The above estimation is based on the dense worst case running times. Since
the matrix is sparse, we can rather assume that the running times for
solving the equation system are about linear with respect to the size~\cite{AmeGLP06,BoGH03,SchWW08}. To estimate the running time for linear programming, we have to resort to
experiments based on the data of Section~\ref{sec:Example1D}. We use the matrices as they would result from starting our algorithm at levels 7 to 14. The results
are shown in Table~\ref{tab:RunningTimes} and
Figure~\ref{fig:RunningTimes}. We used the barrier solver of CPLEX with
additional crossover to recover a basic solution. The running times are
with respect to an Intel Core 2 Quad Core with 2.66 GHz. 
In Table~\ref{tab:RunningTimes}, $\alpha$ seems to be of order $\alpha \approx k$.
Moreover, the results
suggest a growth of the running time that is lower than quadratic. This
seems to be a positive sign, which at least does not rule out a possible
practical effectiveness of our approach. It, however, would require a much
more thorough computational study to reach definite conclusions.

\section*{Conclusion}

As mentioned in the introduction, many issues of the approach presented in
this paper have not yet been resolved and many variations are possible. For
instance, it is obvious that a similar approach could be derived using
other dictionaries, e.g., wavelets, instead of finite element functions.
Furthermore, for practical instances, the solution of the
$\ell_1$-minimization problem becomes an issue. One approach would be to
apply different algorithms, for instance, Orthogonal Matching Pursuit,
see~\cite{DavMZ94,PatRK93b,NeeV07}.  Moreover, the special structure of the
stiffness matrices can be exploited and techniques adapted to the
iterative procedure could be developed.

\section*{Acknowledgment}
We thank O.~Holtz, and J.~Gagelman for many fruitful discussions and
 A. Jensen for help in   the numerical experiments.

\bibliographystyle{mod_siam}

{\small
\bibliography{JokMPY09}
}

\end{document}